\documentclass[twoside]{article}
\usepackage{amsmath,amsfonts,color,enumitem,geometry,latexsym,mathrsfs,upgreek,MnSymbol}
\usepackage[colorlinks,linkcolor=blue,citecolor=red,urlcolor=blue]{hyperref}
\usepackage[amsmath,thmmarks]{ntheorem}
\usepackage[pagewise]{lineno}

\hypersetup{urlcolor=blue, citecolor=red}
\newcommand{\Title}{Topological pressure and equilibrium state for certain correspondences}
\newcommand{\Titleabbr}{\small{Topological pressure and equilibrium state for certain correspondences}}
\date{}
\newcommand{\Keyword}{correspondence; topological pressure; equilibrium state; Ruelle operator; transition probability kernel; bounded measurable potential}
\newcommand{\Fund}{This work is supported by  NSFC (No: 12171400, 12126102, 12426107).}
\newcommand{\Authora}{Yu Zhang}
\newcommand{\AuthorabbrA}{Zhang}
\newcommand{\Authorb}{Yujun Zhu}
\newcommand{\AuthorabbrB}{Zhu}
\newcommand{\Emaila}{exampleA@email.com}
\newcommand{\Emailb}{exampleB@email.com}

\usepackage{graphicx}

\makeatletter
\DeclareFontEncoding{LS2}{}{\noaccents@}
\DeclareFontSubstitution{LS2}{stix}{m}{n}
\DeclareFontFamily{LS2}{stixmathbb}{\skewed@math@base{LS1}{m}{n}}%
\DeclareFontShape{LS2}{stixmathbb}{m}{n} { <-> s*[1.4] stix-mathbb}{}
\makeatother

\usepackage{fancyhdr}
\title{\Title
\footnotetext{
\setlength{\parindent}{2em}
\emph{2020 Mathematics Subject Classification: 37D35, 37D20, 54C60.}\\
\indent \emph{Keywords and phrases:} \Keyword.\\
\indent{\Fund}}}
\author{\Authora \text{ and} \Authorb{}}

\usepackage{enumitem} 
\newtheorem{theorem}{Theorem}[section]

\newtheorem{lemma}[theorem]{Lemma}
\newtheorem{proposition}[theorem]{Proposition}

\newtheorem{remark}[theorem]{Remark}

\theoremstyle{definition}
\newtheorem{definition}[theorem]{Definition}

\newtheorem*{exa*}{Example} 

\newtheorem{thmm}{Theorem}

\theorembodyfont{\normalfont}
\newtheorem*{pf}{Proof}
\newtheorem*{potA}{Proof of Theorem \ref{thm:pressure-m}}
\newtheorem*{potB}{Proof of Theorem \ref{thm:kernel-dense-function}}
\newtheorem*{potC}{Proof of Theorem \ref{thmm-bounded-poten-equi}}

\begin{document}
\vskip 5.1cm
\maketitle
\date{}
\medskip
 {\footnotesize
   \centerline{ School of Mathematical Sciences}
   \centerline{ Xiamen University, Xiamen, 361005, P.R. China}
} 

\bigskip

\begin{abstract}
In \cite{Miller-Akin1999}, Miller and Akin investigated the invariant measures for correspondences, which are also known as upper semi-continuous set-valued maps. Recently, the variational principle and thermodynamic formalism for forward expansive correspondences were studied by Li, Li and Zhang \cite{Xiaoran Li-Zhiqiang Li-Yiwei Zhang2023}. In this paper, the invariant measures and the associated transition probability kernels are explicitly expressed for certain correspondences satisfying the assumptions in \cite{Xiaoran Li-Zhiqiang Li-Yiwei Zhang2023} via the equilibrium states of some particular potentials. Let $T$ be a correspondence on a closed connected Riemannian manifold generated by finite $C^{2}$-expanding endomorphisms. When the generators of $T$ have no coincidence point, a locally H\"{o}lder continuous potential $\phi$ of two variables is defined via the Jacobians of the generators. The pressure of $\phi$  and its equilibrium state  $(\mu, \mathcal{Q})$ are obtained, where $\mu$ is a $T$-invariant measure which is absolutely continuous with respect to the volume and $\mathcal{Q}$ is the associated transition probability kernel satisfying $\mu\mathcal{Q}=\mu$.
For the correspondence $T$ on the torus whose generators have coincidence points, the variational topological pressures for measurable potentials are introduced and the corresponding equilibrium states are considered.
Moreover, the uniqueness of the equilibrium states of  correspondences is considered via the natural extensions.
\end{abstract}

\renewcommand{\sectionmark}[1]{}

\section{Introduction}

Let $X$ be a compact metric space, $f$ a continuous map on $X$ and $\mu$ an $f$-invariant Borel probability measure. The dynamics induced by iterations of a single transformation $f$ are primarily considered from the perspective of measure-theoretic, topological and differentiable  dynamical systems. The measure-theoretic entropy $h_{\mu}(f)$ was presented by
Kolmogorov \cite{Kolmogorov1958} and Sinai \cite{Sinai1959} in 1958, while Adler, Konheim and McAndrew \cite{Adler-Konheim-McAndrew1965} introduced topological entropy $h_{top}(f)$.
Topological pressure $P(f,\varphi)$, a non-trivial weighted topological entropy with respect to a continuous function $\varphi$ on $X$ which is called a \emph{potential} $\varphi$,  is one of the most fundamental dynamical invariants.
 According to the variational principle from \cite{Walters1975},
$P(f,\varphi)= \sup_{\mu}\Big\{h_{\mu}(f)+ \int_{X} \varphi d\mu \Big\}$,
 where $\mu$ ranges over all $f$-invariant measures.
 $f$ may process many invariant measures, some of which may be more interesting than others.
 The \emph{equilibrium state} of $f$ with respect to $\varphi$ is such an invariant measure  that attaches the supremum in the variational principle.

Many classical techniques were developed to find invariant measures, especially for the smooth maps with certain hyperbolicity.
For a $C^{2}$-expanding map $f$ on a closed connected Riemannian manifold $M$, Krzyzewski and Szlenk \cite{Krzyzewski-Szlenk1969} obtained a smooth $f$-invariant measure on $M$ which is absolutely continuous with respect to the volume.
The thermodynamic formalism is a successful approach to find invariant measures with respect to the Jacobian function.
In the 1970s, Sinai \cite{Sinai1972}, Ruelle \cite{Ruelle1978} and Bowen \cite{Bowen1975} considered the thermodynamic formalisms of smooth dynamical systems and  explicitly constructed  equilibrium states for expanding maps and uniformly hyperbolic diffeomorphisms.
In recent years, the thermodynamic formalisms for various systems, such as random dynamical systems, non-uniformly expanding systems, certain non-hyperbolic systems and the correspondences have been extensively investigated, see \cite{Bogenschutz-Gundlach1995}, \cite{P.-D. Liu1998}, \cite{Oliveira2003}, \cite{Yuri1999}, \cite{Xiaoran Li-Zhiqiang Li-Yiwei Zhang2023}, for example.

Correspondences, also known as upper semi-continuous set-valued maps, originated from the study of problems in partial differential equations in the sense of  Hadamard \cite{Hadamard1902}, where solution existence or uniqueness fails under certain data conditions.
 As a natural generalization of single-valued continuous maps, the correspondence has found extensive applications in many fields, including control theory, qualitative physics and complex dynamics.

In recent years, the dynamics of  the correspondences  has attracted more attention.
For example, from the topological and differentiable aspect,
Kelly and Tennant \cite{Kelly-Tennant2017}, as well as Alvin and Kelly \cite{Alvin-Kelly2021}, explored the topological entropies of different types of correspondences.
Pacifico and Vieitez  \cite{Pacifico-Vieitez2017} considered the expansiveness, Lyapunov exponents and entropy for correspondences.
Meanwhile, Kennedy and Nall \cite{Kennedy-Nall2018} as well as Ingram and Mahavier \cite{Ingram-Mahavier2006} studied  correspondences via their inverse limits. Entropy formulas and entropy inequalities for the correspondences generated by discrete group actions on $M$ were also investigated (see  \cite{Zhu2017}, \cite{Zhu-Zhang2022}, \cite{Wang-zhang-Zhu2023} and \cite{Wang-Zhang-Zhu2025}, etc).
From the measure-theoretic aspect,
Aubin, Frankowska and Lasota \cite{Aubin-Frankowska-Lasota1991} and Miller \cite{{Miller1995}} introduced two concepts of invariant measures respectively and gave the corresponding existence theorems for a correspondence on a compact metric space $X$.
Miller and Akin \cite{Miller-Akin1999} introduced two more concepts of invariant measures and showed that all  these four definitions are equivalent. Among these definitions, the one with the simplest representation is as follows. Let $T$  be a  correspondence on $X$ and $\mu$ a Borel probability measure of $X$. $\mu$ is said to be $T$-invariant if $\mu (A) \leq \mu (T^{-1}(A))$, $\forall A \in \mathcal{B}(X)$. Clearly, it is impossible to introduce a meaningful concept of measure-theoretic entropy of $T$ only using $\mu$. \cite{Miller-Akin1999} showed that $\mu$ is $T$-invariant if and only if there exists a transition probability kernel $\mathcal{Q}$ supported by $T$ such that $\mu\mathcal{Q}=\mu$.
Recently, Li, Li and Zhang \cite{Xiaoran Li-Zhiqiang Li-Yiwei Zhang2023}
introduced a kind of measure-theoretic entropy for $T$ using $\mu$ and $\mathcal{Q}$ together. Moreover, when $T$ is forward expansive, a variational principle for the pressure of a potential of two variables was  established  and the existence and uniqueness of the related equilibrium state  $(\mu, \mathcal{Q})$ was considered in \cite{Xiaoran Li-Zhiqiang Li-Yiwei Zhang2023}.

The main purpose of this paper is to give the explicit expression of the invariant measures and the associated transition probability kernels for certain correspondences satisfying the assumptions in \cite{Xiaoran Li-Zhiqiang Li-Yiwei Zhang2023} via the equilibrium states of some appropriate potentials.  Let $M$ be a closed connected Riemannian manifold and $f_1,f_{2}, \cdots, f_k$ be $C^{2}$-expanding endomorphisms. A correspondence $T$ can be defined by $T(x)=\{f_{j}(x):1\le j\le k\}$.
When each pair of the generators has no coincidence point, we will focus on the following potential,
 \begin{equation}\label{volume-m}
  \phi(x_{1},x_{2})=-\log \text{Jac}(f_{j}(x_{1})),\;\;\text{ if }x_{2}=f_{j}(x_{1}), \;\;1 \leq j \leq k,
 \end{equation}
where  $\text{Jac}(f(x))$ is the Jacobian of the linear map $Df(x)\colon T_{x}M \to T_{f(x)}M$. If $k=1$, then the topological pressure of $f$ with respect to such $\phi$ equals 0, the equilibrium state is  absolutely continuous with respect to volume measure, and it is just the unique Gibbs (SRB) measure.
We will study  the topological pressures and  equilibrium states of the correspondence $T$ with respect to such $\phi$, and then explore the uniqueness of the equilibrium state via the corresponding natural extension (i.e., the left shift map on the orbit space of the correspondence).

We first show that if a correspondence $T$ is generated by finite $C^{2}$-expanding endomorphisms with no coincidence point, then $T$ is distance-expanding and the potential $\phi$ is locally H\"{o}lder.
Inspired by the work of Hurley \cite{Hurley1995} and Nitecki and Przytycki \cite{Nitecki-Przytycki1999}, we then analyze the local contracting behavior along backward orbits. We conclude that the topological pressure $P(T,\phi)$ equals $\log k$, where $k$ is the number of generators of $T$ (Theorem \ref{thm:pressure-m}). The proof's key is to show that $P(T,\phi)$ can be calculated via the exponential growth rate of the number of backward $n$-orbits weighted by the Jacobian at an arbitrary point $x \in M$.

Combining the fact that $P(T,\phi)=\log k$ with the variational principle for correspondences established by Li, Li and Zhang \cite{Xiaoran Li-Zhiqiang Li-Yiwei Zhang2023}, we obtain the equilibrium state $(\mu,\mathcal{Q})$ for a certain correspondence $T$ on $M$ (Theorem \ref{thm:kernel-dense-function}). Here, $\mu$ is a $T$-invariant smooth probability measure, and $\mathcal{Q}=\frac{1}{k}\sum_{j=1}^{k}\delta_{f_{j}}$.
It is well known that for a $C^{2}$-expanding endomorphism on $M$, the thermodynamic formalism is crucial for obtaining the density function and the smooth invariant measure. However, for the correspondence $T$ generated by finite $C^{2}$-expanding endomorphisms, the thermodynamic formalism for a single-valued map is inapplicable.
To address this issue, we construct a type of Ruelle operator to get the density function for the correspondence. The density function is the eigenfunction of the Ruelle operator with eigenvalue $k$.
Thus, we obtain the smooth probability measure $\mu$ on $M$ and show that it is $\mathcal{Q}$-invariant, where $\mathcal{Q}=\frac{1}{k}\sum_{j=1}^{k}\delta_{f_{j}}$.
Next, we connect the measure-theoretic entropy of the transition probability kernel $\mathcal{Q}$ for the smooth probability measure $\mu$ to the corresponding skew product transformation. The fiber-entropy of the skew product can be calculated using the Jacobian of each generator $f_{j}$, $1\leq j\leq k$, which makes $(\mu, \mathcal{Q})$ an equilibrium state of $P(T,\phi)$.
We remark in Section 6 that the equilibrium state obtained in Theorem \ref{thm:kernel-dense-function} is unique in an appropriate  sense.

For the correspondence $T$ which is generated by finite $C^{2}$-expanding endomorphisms with coincidence points, things will become more complicated. We consider a special case.
Let $A_{j}, 1\le j\le k$ be all non-singular $m\times m$ integer matrices and  all eigenvalues of each $A_{j}$ be of modulus greater than 1. Then there is an induced correspondence $T$ on torus  $\mathbb{T}^{m}$.  In this case, we are not able to define a continuous function as in (\ref{volume-m}). However, we can modify the definition to get a measurable function $\phi$ as in (\ref{spe-poten}).  So, we turn to consider  the variational topological pressure $P^{v}(T, \phi)$ for a real-valued bounded measurable function. We will show that $P^{v}(T, \phi) = \log k$ and $(\text{Leb},\frac{1}{k}\sum_{j=1}^{k}\delta_{f_{j}})$ is an equilibrium state for the variational topological pressure (Theorem \ref{thmm-bounded-poten-equi}). By analyzing the local contracting behavior along backward orbits to calculate the number of elements in a spanning set, we derive the inequality $P^{v}(T, \phi) \leq \log k$. Note that
for each pair of generators $f_1,f_{2}, \cdots, f_k$, the set of their coincidence points is of zero Lebesgue measure. Then, using the same idea as in the proof of Theorem \ref{thm:kernel-dense-function}, we obtain the opposite inequality $P^{v}(T, \phi)\geq\log k$.


The paper is organized as follows. In Section 2, we introduce the basic notions for correspondences and state the main results. The proofs of Theorem \ref{thm:pressure-m}, Theorem \ref{thm:kernel-dense-function} and Theorem \ref{thmm-bounded-poten-equi} are given in Section 3,  Section 4 and Section 5, respectively.
Section 6 is devoted to discussing the uniqueness of equilibrium states for correspondences via the corresponding  natural extension systems.


\section{Notations and statement of main results}


\subsection{Correspondences}\label{subsection2-1}
Throughout this paper we denote the set of integers, the set of positive integers and the set of non-negative integers by $\mathbb{Z},\mathbb{N}$ and $\mathbb{Z}_+$ respectively. Let $\widehat{\mathbb{N}}:=\mathbb{N}\cup \{\omega\}$. Here $\omega$ is the least infinite ordinal.
Let $(X, d)$ be a compact metric space.
For each $n \in \mathbb{N}$, equip $X^{n}:=\{(x_{1},x_{2},\cdots, x_{n}):x_{i}\in X \text{ for all } 1 \leq i \leq n\}$
with the metric $d_{n}$ given by
$$
d_{n}((x_{1},x_{2},\cdots, x_{n}), (y_{1},y_{2},\cdots, y_{n})):=\max_{1 \leq i \leq n}d(x_{i},y_{i})
$$for $(x_{1},x_{2},\cdots, x_{n}), (y_{1},y_{2},\cdots, y_{n}) \in X^{n}$.
Similarly, equip $X^{\omega}:=\{(x_{1},x_{2},\cdots):x_{i}\in X \text{ for all } i \in \mathbb{N}\}$
with the metric $d_{\omega}$ given by
\begin{equation*}\label{orbit metric}
d_{\omega}((x_{1},x_{2},\cdots), (y_{1},y_{2},\cdots)):=\sum_{i=1}^{+\infty}\frac{d(x_i,y_i)}{2^{i}(d(x_i,y_i)+1)}
\end{equation*}for  $(x_{1},x_{2},\cdots), (y_{1},y_{2},\cdots) \in X^{\omega}$.
With the metrics $d_{n}$ for $n \in \widehat{\mathbb{N}}$, the topologies of $X^{n}$ induced by these metrics are the product topologies.

Let $\sigma\colon X^{\omega}\to X^{\omega}$ be the standard left shift map defined by
 $$
 \sigma(x_{1},x_{2},x_{3},\cdots)=(x_{2},x_{3},\cdots)
$$for  $(x_{1},x_{2},x_{3}\cdots) \in X^{\omega}$.
Clearly, $\sigma$ is a continuous surjection on $X^{\omega}$.

We define $2^X$ to be the set of all non-empty, closed
subsets of $X$. Let $T\colon X\to 2^X$ be a set-valued map with $\bigcup_{x\in X}T(x)=X$, then its inverse $T^{-1}$ is a set-valued map defined as $T^{-1}(x)=\{y\colon x\in T(y)\}, x\in X$. Let $T'\colon X\to 2^X$ be also a set-valued map, then the composition of $T$ and $T'$ is a set-valued map defined by $T'\circ T(x)=\cup_{y\in T(x)}T'(y), x\in X$. For $n\in \mathbb{Z}_+$, denote $T^{n}$ the $n$-th composition of $T$.
We say that a set-valued map $T\colon X\to 2^X$ is \emph{upper semi-continuous} if for every
open set $U \subset X$, the set $\{x\in X : T(x) \subset U\}$ is an open subset of $X$. Once $T$ is upper semi-continuous then so is $T^{n}$ for each $n\in \mathbb{Z}_+$.

In this paper, we call an upper semi-continuous set-valued map $T\colon X\to 2^X$ a \emph{correspondence} on $X$.
For a correspondence $T\colon X\to 2^X$, we define the set of \emph{forward n-orbit space} as
$$
\mathcal{O}_{n}{(T)}:=\{\underline{x}=(x_{1},x_{2},\cdots, x_{n})\in X^{n}:x_{i+1}\in T(x_{i}) \text{ for each } i \in \{1,2,\cdots,n-1\}\} \subset X^{n}.
$$
The \emph{forward orbit space} $\mathcal{O}_{\omega}{(T)}$ induced by $T$ is given by
\begin{equation*}\label{orbit space}
\mathcal{O}_{\omega}{(T)}:=\{\underline{x}=(x_{1},x_{2},\cdots)\in X^{\omega}:x_{i+1}\in T(x_{i}) \text{ for each } i \in \mathbb{N}\} \subset X^{\omega}.
\end{equation*}
A well-known result is that a set-valued map $T$ is upper semi-continuous if and only if $\mathcal{O}_{2}{(T)}=\{(x_{1}, x_{2}) \in X^{2}: x_{2}\in T(x_{1})\}$ is closed in $X^{2}$ (see \cite{Xiaoran Li-Zhiqiang Li-Yiwei Zhang2023}, for example).
Clearly, $\mathcal{O}_{\omega}{(T)}$ is $\sigma$-invariant. Let $\sigma_T$ be the restriction of $\sigma$ to $\mathcal{O}_{\omega}{(T)}$. We call $(\mathcal{O}_{\omega}{(T)}, \sigma_T)$  the \emph{natural extension} of $(X, T)$.

Let $\phi\colon \mathcal{O}_{2}{(T)}\to \mathbb{R}$ be a continuous function on $(\mathcal{O}_{2}{(T)}, d_{2})$. Denote by $\widetilde{\phi}\colon \mathcal{O}_{\omega}{(T)}\to \mathbb{R}$ the function given by
\begin{equation}\label{local constant}
\widetilde{\phi}(x_{1},x_{2},x_{3}\cdots):=\phi(x_{1},x_{2}).
\end{equation}Clearly, $\widetilde{\phi}$ is also a continuous function on $(\mathcal{O}_{\omega}{(T)}, d_{\omega})$.

Let $M$ be an m-dimensional  closed connected Riemannian manifold without boundary. Then $M$ has a Riemannian structure, which induces a volume measure $\mathbf{m}$ on $M$.
In this paper, we will consider a particular case of correspondences which is generated by finite continuous or smooth surjective single-valued maps.

\begin{definition}\label{finite-generated}
 Let $f_1,f_2, \cdots, f_k$ be continuous (resp. smooth) maps on $X$ (resp. $M$). We can define a correspondence $T\colon X\to 2^X$ as follows,
\begin{equation*}
T(x)\colon =\{y\colon y=f_j(x) \text{ for some } 1\le j\le k\},
\end{equation*}which is called \emph{finitely-generated correspondence with the generators $f_1,f_2, \cdots, f_k$}.
\end{definition}
Clearly, $\mathcal{O}_{2}{(T)}=\bigcup_{j=1}^{k}\mathcal{O}_{2}{(f_{j})}$ is closed in $X^2$ under the product topology.


\subsection{Topological pressure for correspondences}
Now, we recall some  definitions and properties of pressure for correspondences.
For $\varepsilon >0$, and a compact metric space $(Y,\rho)$,  a subset $E \subset Y$ is said to be an $(\rho,\varepsilon)$-separated set of $(Y,\rho)$, if any $x,y \in E$, $x\neq y$ implies $\rho(x,y)> \varepsilon$.
A subset $F \subset Y$ is said to be an $(\rho,\varepsilon)$-spanning set of $(Y,\rho)$, if for any $x \in Y$, there exists $y \in F$ such that
$\rho(x,y)\leq \varepsilon$.

Let $f$ be a  continuous map on a compact metric space $(X,d)$ and $\varphi\colon X \to \mathbb{R}$ be a continuous function. For any $n\in \mathbb{N}$, define a
new metric $d_{n}^{\prime}$ on $X$ by
$$
d_{n}^{\prime}(x,y):=\max_{0 \leq i \leq n-1} d(f^{i}(x),f^{i}(y)).
$$

The\emph{ topological pressure} of $f$ with respect to the continuous potential $\varphi \colon X \to \mathbb{R}$ is defined by
\begin{align}\label{single-preuusre}
 P(f,\varphi)=\lim\limits_{\varepsilon\rightarrow 0} \limsup_{n\rightarrow +\infty}\frac{1}{n}\log \inf_{F} \big\{\sum_{x\in F}\exp(\sum_{i=0}^{n-1}\varphi(f^{i}x))\big\},
\end{align}
where the  infimum is taken over all $(d_{n}^{\prime},\varepsilon)$-spanning sets of $(X,d)$. If 0 denotes the constant function with valued 0, then $ P(f,0)=h_{top}(f)$, the \emph{topological entropy} of $f$ (see \cite{Walters1982}).

Let $T$ be a correspondence on a compact metric space $(X,d)$ and $\phi \colon \mathcal{O}_{2}(T) \to \mathbb{R}$ be a continuous function.
For each $n \in \mathbb{N}$, define $S_{n}\phi \colon \mathcal{O}_{n+1}(T) \to \mathbb{R}$ is given by
\begin{eqnarray*}\label{summ}
  S_{n}\phi(x_{1},x_{2},\cdots,x_{n+1})=\sum_{i=1}^{n}\phi(x_{i},x_{i+1}).
\end{eqnarray*}
For each $n \in \mathbb{N}$, write
$$
S_{n}(T,\phi,\varepsilon)=\sup\big\{\sum_{\underline{x}\in E}\exp(S_{n}\phi(\underline{x}))|{E \text{ is a }(d_{n+1},\varepsilon)\text{-separated set of }(\mathcal{O}_{n+1}(T),d_{n+1})}\big\},
$$and
$$
R_{n}(T,\phi,\varepsilon)=\inf\big\{\sum_{\underline{x}\in F}\exp(S_{n}\phi(\underline{x}))|{F \text{ is a }(d_{n+1},\varepsilon)\text{-spanning set of }(\mathcal{O}_{n+1}(T),d_{n+1})}\big\}.
$$
\begin{definition}(Definition $4.6$ in  \cite{Xiaoran Li-Zhiqiang Li-Yiwei Zhang2023})
 Let $T$ be a correspondence on a compact metric space $(X,d)$ and $\phi \colon \mathcal{O}_{2}(T) \to \mathbb{R}$ be a continuous function. Then the \emph{topological pressure} of $T$ with respect to the potential $\phi$ is defined by
\begin{eqnarray}\label{toppre}
   P(T,\phi)=\lim\limits_{\varepsilon\rightarrow 0} \limsup_{n\rightarrow+\infty}\frac{1}{n}\log S_{n}(T,\phi,\varepsilon).
\end{eqnarray}
\end{definition}
Moreover, $P(T,\phi)=\lim\limits_{\varepsilon\rightarrow 0} \limsup_{n\rightarrow+\infty}\frac{1}{n}\log R_{n}(T,\phi,\varepsilon)$ (see \cite{Xiaoran Li-Zhiqiang Li-Yiwei Zhang2023}). Note that this definition of topological pressure  naturally generalizes the definition of topological pressure of a single-valued continuous map as in $(\ref{single-preuusre})$.

\begin{proposition}(Theorem $4.9$ in \cite{Xiaoran Li-Zhiqiang Li-Yiwei Zhang2023})
Let $T\colon X\to 2^X$ be a correspondence on $(X, d)$ and $\sigma_T\colon \mathcal{O}_{\omega}{(T)}\to \mathcal{O}_{\omega}{(T)}$ be its natural extension. Let $\phi \colon \mathcal{O}_{2}(T) \to \mathbb{R}$ be a continuous function.
 Then
\begin{equation}\label{topological pressure}
 P(T,\phi)=P(\sigma_T,\widetilde{\phi}),
\end{equation}where $\widetilde{\phi}\colon \mathcal{O}_{\omega}{(T)}\to \mathbb{R}$ is a continuous function on $(\mathcal{O}_{\omega}{(T)},d_{\omega})$ with $\widetilde{\phi}(x_{1},x_{2},x_{3},\cdots):=\phi(x_{1},x_{2})$.
\end{proposition}

In particular, if $\phi=0$, then $P(T,0)=P(\sigma_T,0)=h_{top}(\sigma_T)$. It follows that when $\phi$ vanishes, the notion of topological pressure  coincides with the notion of topological entropy of correspondences (see \cite{Kelly-Tennant2017}, for example).

\begin{definition}\label{expansive}
A correspondence $T\colon X\to 2^X$ is \emph{forward expansive} if there is $c >0$ (call \emph{expansive constant}) such that for each pair of distinct orbits $(x_{1},x_{2},\cdots), (y_{1},y_{2},\cdots) \in \mathcal{O}_{\omega}{(T)}$, we have $d(x_{i},y_{i})>c$ for some $i \in \mathbb{N}$.
\end{definition}

We denote $|\cdot|$
and $d(\cdot,\cdot)$, respectively, the norm on the tangent bundle $TM$ and the metric on $M$ induced by the Riemmanian
metric (except mentioned otherwise). Let $\varepsilon>0$. We denote $B(x,\varepsilon)$ the open $\varepsilon$ ball about $x$ in $M$.

\begin{definition}[\cite{Qian-Xie-Zhu2009}]\label{def:hp}
A map $f\colon M \to M$ is called an \emph{expanding map} on $M$, if there exist $\lambda_{0}>1$ and $c_{0}> 0$ such that for all $x \in M$ and $n \geq 0$,
  \begin{equation*}
  |Df^n(x)v|\geq c_{0}\lambda_{0}^{n}|v|, \text{ for } v \in T_{x}M.
 \end{equation*}
\end{definition}

\begin{remark}
Via a standard change of Riemmanian metric we may, and will always, assume that $c_{0}=1$ (see \cite{Wen2016} for the invertible case).
\end{remark}

Throughout this paper, we assume $M$ is oriented. For any $C^{2}$-expanding map $f$ on $M$, let $\text{deg}(f)$ denote the degree of $f$. Note that $\text{deg}(f)$ equals the number of preimages of any point in $M$. Recall that  for any $x \in M$, $\text{Jac}(f(x))$ denotes the Jacobian of the linear map $Df\colon T_{x}M \to T_{f(x)}M$.

We now consider the forward expansive correspondence. Let $f$ and $g$ be continuous maps on $X$, $x$ is called the \emph{coincidence} of $f$ and $g$ if $f(x)=g(x)$. Denote the set of the coincidence points of $f$ and $g$ by $\text{Coinc}(f,g)$.

\begin{thmm}\label{thm:pressure-m}
Let $T$ be a finitely-generated correspondence with finite $C^{2}$-expanding endomorphisms $f_{1},f_{2}$, $\cdots, f_{k}$ on $M$. Assume $\text{Coinc}(f_{j},f_{j^{\prime}})=\emptyset$ for $1 \leq j\neq j^{\prime} \leq k$. Let $\phi\colon \mathcal{O}_{2}{(T)}\to \mathbb{R}$ be the potential as in (\ref{volume-m}). Then we have
 $$P(T,\phi)=\log k.$$
\end{thmm}

\subsection{Transition probability kernels and equilibrium states for correspondences}

In this subsection, we recall the definition and some properties of transition probability kernels, which are also called Markovian transition kernels (see \cite{Meyn-Tweedie2012} and \cite{Revuz1984}, for example). In fact, a transition probability kernel on a  measurable space $X$ assigns each $x \in X$ a probability measure on $X$ (Definition \ref{tra-pro-kernel}).
 And we also recall how a transition probability kernel pushes forward a probability measure (Definition \ref{pushforward probability measure}) and how  the
  transition probability kernel with a probability measure generates the Markov process (Definition \ref{Markov process}).
The following three  definitions are from  Definitions $5.1$, $5.3$, and $5.6$ in \cite{Xiaoran Li-Zhiqiang Li-Yiwei Zhang2023}.

\begin{definition}\label{tra-pro-kernel}
  Let $(X,d)$ be a compact metric space and $\mathcal{B}(X)$ be the Borel $\sigma$-algebra on $X$. A \emph{transition probability kernel} on $X$ is a map $\mathcal{Q}\colon X\times\mathcal{B}(X) \to (0,1]$ satisfying the following two properties:

  $(1)$ For every $x \in X$, the map $\mathcal{B}(X) \ni A \longmapsto \mathcal{Q}(x,A):=\mathcal{Q}_{x}(A)$ is a Borel probability measure on $(X, \mathcal{B}(X))$.

   $(2)$ For every $A \in \mathcal{B}(X)$, the map $X \ni x \longmapsto \mathcal{Q}_{x}(A)$ is Borel measurable.
\end{definition}

\begin{definition}\label{support}
  Let $T$ be a correspondence on $(X,d)$ and $\mathcal{Q}$ be a transition probability kernel on $(X,\mathcal{B}(X))$. We say $\mathcal{Q}$ is \emph{supported} by $T$ if the measure $\mathcal{Q}_{x}$ is supported on the closed set $T(x)$ (i.e, $\mathcal{Q}_{x}(T(x))=1$) for every $x \in X$.
\end{definition}

\begin{definition}\label{pushforward probability measure}
  Let $(X,\mathcal{B}(X))$ be a Borel measure space, $\mu \in \mathcal{M}(X)$ be a Borel probability measure on $X$ and $\mathcal{Q}$ be a transition probability kernel on $X$. The \emph{pushforward probability measure} $\mu\mathcal{Q}$ on $X$ of $\mu$ by $\mathcal{Q}$ is given by
  $$
  (\mu\mathcal{Q})(A):=\int_{X} \mathcal{Q}(x,A) d\mu(x),\;\;\; \forall A \in \mathcal{B}(X).
  $$We say that a Borel probability measure $\mu$ on $X$ is \emph{$\mathcal{Q}$-invariant} if $\mu\mathcal{Q}=\mu$. Denote by $\mathcal{M}(X, \mathcal{Q})$ the set of all $\mathcal{Q}$-invariant Borel probability measures on $X$. Denote by $\mathcal{M}_{\sigma_{T}}(\mathcal{O}_{\omega}(T))$ the set of all $\sigma_{T}$-invariant Borel probability measures on $\mathcal{O}_{\omega}(T)$.

\end{definition}

\begin{theorem}(\cite{Miller-Akin1999})
  For a  correspondence $T\colon X\to 2^X$ on a compact metric space $(X,d)$ and $\mu \in \mathcal{M}(X)$,  the following  are equivalent definitions of $\mu$ being an invariant measure for  $T$:

  $(1)$ $\mu (A) \leq \mu (T^{-1}(A))$, $\forall A \in \mathcal{B}(X)$.

  $(2)$ There exists a transition probability kernel $\mathcal{Q}$ supported by $T$ such that $\mu\mathcal{Q}=\mu$.

  $(3)$ There exists a measure $\mu_{12} \in \mathcal{M}(X \times X)$, supported on $\mathcal{O}_{2}(T)$, such that $\mu_{12}\circ {\widetilde{\pi}_{1}}^{-1}=\mu_{12}\circ {\widetilde{\pi}_{2}}^{-1}=\mu$.

  $(4)$ There exists a  probability measure $\nu \in \mathcal{M}_{\sigma_{T}}(\mathcal{O}_{\omega}(T))$  such that $\nu \circ {\widetilde{\pi}_{1}}^{-1}=\mu$.
 \end{theorem}
Therefore, if $\mu \in \mathcal{M}(X,\mathcal{Q})$ and $\mathcal{Q}$ is supported by $T$, then $\mu$ is an invariant measure for $T$.

\begin{remark}\label{Relation-invariant-Markov}
If $\mu$ is an invariant measure for $T$, then there exist measure $\mu_{j} \in \mathcal{M}(X)$, $1\leq j\leq k$, and a  probability vector $\{\alpha_{1},\alpha_{2}, \cdots, \alpha_{k}\}$, each $\alpha_{j}>0$, such that $\mu=\sum_{j=1}^{k}\alpha_{j}\mu_{j}=\sum_{j=1}^{k}\alpha_{j}(\mu_{j}\circ f_{j}^{-1})$.
In fact, if $\mu=\mu_{12}\circ {\widetilde{\pi}_{1}}^{-1}=\mu_{12}\circ {\widetilde{\pi}_{2}}^{-1}$, take $\alpha_{j}=\mu_{12}(\mathcal{O}_{2}(f_{j}))$, and $\mu_{j}=\frac{1}{\alpha_{j}}(\mu_{12}|_{\mathcal{O}_{2}(f_{j})})\circ {\widetilde{\pi}_{1}}^{-1}$, where $\mu_{12}|_{\mathcal{O}_{2}(f_{j})}(A)=\mu_{12}(A \cap \mathcal{O}_{2}(f_{j})), \forall A \in \mathcal{B}(X\times X)$. On one hand, $\mu=\mu_{12}\circ {\widetilde{\pi}_{1}}^{-1}=\sum_{j=1}^{k}\alpha_{j}(\frac{1}{\alpha_{j}}(\mu_{12}|_{\mathcal{O}_{2}(f_{j})})\circ {\widetilde{\pi}_{1}}^{-1})=\sum_{j=1}^{k}\alpha_{j}\mu_{j}$. On the other hand, $\mu=\mu_{12}\circ {\widetilde{\pi}_{2}}^{-1}=\sum_{j=1}^{k}\alpha_{j}(\frac{1}{\alpha_{j}}(\mu_{12}|_{\mathcal{O}_{2}(f_{j})})\circ {\widetilde{\pi}_{2}}^{-1})=\sum_{j=1}^{k}\alpha_{j}(\mu_{j}\circ f_{j}^{-1})$.
\end{remark}

For the finitely-generated correspondence $T$ with finite continuous maps $f_{1},f_{2}, \cdots, f_{k}$ on $X$,
$\mathcal{Q}_{x}=\sum_{j=1}^{k}P_{f_{j}}(x)\delta_{f_{j}(x)}$, $\forall x \in X$, where $\delta_{f_{j}(x)}$ refers to the Dirac measure on $X$ at the point $f_{j}(x)$, and each $P_{f_{j}}$ is a Borel measurable function on $X$ with $\sum_{j=1}^{k}P_{f_{j}}(x)=1$, $0 \leq P_{f_{j}}(x)\leq 1$, $1 \leq j \leq k$, $\forall x \in X$, represent all of the transition probability kernels supported by $T$ on $X$.
For simplicity, we denote $P_{x,f_{j}}:=P_{f_{j}}(x)$, $\forall x \in X$, $1 \leq j \leq k$.
If $\mu \mathcal{Q}=\mu$, then $\alpha_{j}=\int_{X}P_{f_{j}}(x) d\mu(x)$, and $d\mu_{j}(x)=\frac{1}{\alpha_{j}}P_{f_{j}}(x)d\mu(x), \forall x \in X$, $1\leq j \leq k$. Hence, $\mu=\sum_{j=1}^{k}\alpha_{j}\mu_{j}=\sum_{j=1}^{k}\alpha_{j}(\mu_{j}\circ f_{j}^{-1})$.

Denote $\mathcal{B}(X^{\omega}):=\mathcal{B}(\prod_{n \in \mathbb{Z}_{+}}X)$ the Borel $\sigma$-algebra on $X^{\omega}$ generated by $\bigcup_{n=0}^{+\infty}\{X^{n}\times A \times X^{\omega}: A \in \mathcal{B}(X)\}$. For each $A_{n+1} \subset X^{n+1}$  and each $(x_{1}, x_{2}, \cdots, x_{n})\in X^{n}$, write
\begin{equation*}\label{projection-n}
  \pi_{n+1}(x_{1},x_{2},\cdots, x_{n};\;A_{n+1})=\{x_{n+1} \in X; \;(x_{1}, x_{2}, \cdots, x_{n}, x_{n+1}) \in A_{n+1}\}.
\end{equation*}

\begin{definition}(Definition $5.14$, $5.18$ in \cite{Xiaoran Li-Zhiqiang Li-Yiwei Zhang2023})\label{Markov process}
  Let $\mathcal{Q}$ be a transition probability kernel on $(X,\mathcal{B}(X))$.

  $(1)$
  Define the transition probability kernel $\mathcal{Q}^{[n]}$ from $X$ to $X^{n+1}$ inductively on $n \in \mathbb{Z}_{+}$ as follows: for each $x \in X$, $\mathcal{Q}_{x}^{[0]}=\delta_{x}$. If $\mathcal{Q}^{[n-1]}$ has defined for some $n \in \mathbb{N}$, we define $\mathcal{Q}^{[n]}$ as
  \begin{equation*}\label{stochastic process}
    \mathcal{Q}^{[n]}(x,A_{n+1}):=\int_{X^{n}} \mathcal{Q}(x_{n}, \pi_{n+1}(x_{1},x_{2},\cdots, x_{n};\;A_{n+1}))  d\mathcal{Q}_{x}^{[n-1]}(x_{1},x_{2},\cdots,x_{n}),
  \end{equation*}where $A_{n+1} \in \mathcal{B}(X^{n+1})$.

  $(2)$ Define the transition probability kernel $\mathcal{Q}^{\omega}$ from $X$ to $X^{\omega}$ as the unique transition probability kernel from $X$ to $X^{\omega}$ with the property that for each $x \in X$, each $n \in \mathbb{Z}_{+}$, each $A \in \mathcal{B}(X^{n+1})$,
  $$
  \mathcal{Q}^{\omega}(x, A\times X^{\omega})=\mathcal{Q}^{[n]}(x,A).
  $$
\end{definition}

 \begin{remark}\label{forward-infite}
   $(1)$
   For each $\mu \in \mathcal{M}(X)$, Definition  \ref{pushforward probability measure} and \ref{Markov process} imply
   $$
   (\mu\mathcal{Q}^{\omega})(A \times X^{\omega})=(\mu\mathcal{Q}^{[n]})(A), \;\;\; \forall A \in \mathcal{B}(X^{n+1}),\; \forall n \in \mathbb{Z}_{+}.
   $$

$(2)$
 If $\mu \in \mathcal{M}(X,\mathcal{Q})$ and $\mathcal{Q}$ is supported by $T$, then $\mu\mathcal{Q}^{\omega}$ is
 $\sigma_{T}$-invariant (see \cite{Xiaoran Li-Zhiqiang Li-Yiwei Zhang2023}, for example).
 \end{remark}

Li, Li and Zhang \cite{Xiaoran Li-Zhiqiang Li-Yiwei Zhang2023} introduced a kind of measure-theoretic entropy $h_{\mu}(\mathcal{Q})$ for the transition probability kernel $\mathcal{Q}$ and the $\mathcal{Q}$-invariant Borel probability measure $\mu$, which generalizes the classical measure-theoretic entropy of single-valued maps.

\begin{definition}\label{parti-entropy and entropy}(Definition $5.21$, $5.22$ in \cite{Xiaoran Li-Zhiqiang Li-Yiwei Zhang2023})
  Let $\mathcal{Q}$ be a transition probability kernel on $(X,\mathcal{B}(X))$, $\mu$ be a $\mathcal{Q}$-invariant Borel probability measure on $X$, and $\mathcal{A}$ be a finite measurable partition of $X$.  Then $h_{\mu}(\mathcal{Q},\mathcal{A})$, the measure-theoretic entropy of $\mathcal{Q}$ with respect
 to the partition $\mathcal{A}$, is defined as
\begin{equation*}
h_{\mu}(\mathcal{Q},\mathcal{A}):=\lim_{n\rightarrow +\infty}\frac{1}{n}H_{\mu \mathcal{Q}^{[n-1]}}(\mathcal{A}^{n}),
\end{equation*}where $\mathcal{A}^{n}=\{A_{1}\times \cdots \times A_{n}: A_{j} \in \mathcal{A}, 1 \leq j \leq n\}$, and $H_{\mu}(\mathcal{A})=-\sum_{A \in \mathcal{A}} \mu(A)\log\mu(A)$. Then the measure-theoretic entropy $h_{\mu}(\mathcal{Q})$ (of $\mathcal{Q}$ for $\mu$) is given by
\begin{equation*}
  h_{\mu}(\mathcal{Q}):=\sup_{\mathcal{A}}h_{\mu}(\mathcal{Q},\mathcal{A}),
\end{equation*}where $\mathcal{A}$ ranges over all finite measurable partitions of $X$.
\end{definition}

\begin{remark}
  From Remark \ref{forward-infite} $(2)$, the \emph{measure-theoretic entropy} $h_{\mu\mathcal{Q}^{\omega}}(\sigma_{T})$ of $\sigma_{T}$ with respect to $\mu\mathcal{Q}^{\omega}$ is well-defined (see \cite{Walters1982}).
  If $\mathcal{Q}$ is supported by a correspondence $T$ and $\mu \mathcal{Q}=\mu$, then
  $h_{\mu}(\mathcal{Q})=h_{\mu\mathcal{Q}^{\omega}}(\sigma_{T})$ (see Theorem $5.24$ of \cite{Xiaoran Li-Zhiqiang Li-Yiwei Zhang2023}).
\end{remark}

In \cite{Xiaoran Li-Zhiqiang Li-Yiwei Zhang2023},  Li, Li and Zhang  used transition probability kernels to  establish a complete  variational principle and prove the existence of the related equilibrium states for forward expansive correspondences as follows:

\begin{theorem}\label{variational principle-from}(Part $(i)$ of Theorem A in \cite{Xiaoran Li-Zhiqiang Li-Yiwei Zhang2023})
  Let $(X,d)$ be a compact metric space, $T$ be a forward expansive correspondence on $X$, and $\phi\colon \mathcal{O}_{2}(T)\to \mathbb{R}$ be a continuous function. Then the Variational Principle holds:
  \begin{equation}\label{variational principle-correspondence}
    P(T,\phi)= \sup_{(\mu,\mathcal{Q})}\bigg\{h_{\mu}(\mathcal{\mathcal{Q}})+ \int_{X}\int_{T(x_{1})} \phi(x_{1},x_{2})d\mathcal{Q}_{x_{1}}(x_{2})d\mu(x_{1})\bigg\} \in \mathbb{R},
  \end{equation}where $\mathcal{Q}$ ranges over all transition probability kernels on $X$ supported by $T$ and $\mu$
ranges over all $\mathcal{Q}$-invariant Borel probability measures on $X$.
\end{theorem}

\begin{remark}
  If a transition probability kernel $\mathcal{Q}$ on $X$ supported by $T$ and a $\mathcal{Q}$-invariant Borel probability measure $\mu$
on $X$ satisfy the equality $(\ref{variational principle-correspondence})$, then the pair $(\mu, \mathcal{Q})$ is called the equilibrium state for the correspondence $T$ and the potential $\phi$. According to part $(ii)$ of Theorem $A$ in \cite{Xiaoran Li-Zhiqiang Li-Yiwei Zhang2023}, an equilibrium state $(\mu, \mathcal{Q})$ exists for $P(T,\phi)$.
 \end{remark}

For any $n \in \mathbb{N}$, $x \in M$, let
\begin{align*}
 \mathcal{O}_{n+1}(x)=&\{\underline{x}=(x_{1},\cdots, x_{n}, x_{n+1}) \in \mathcal{O}_{n+1}(T):
 x_{n+1}=x\},
\end{align*}and
$$\Phi_{n}(x)=\sum_{\underline{x}=(x_{1},\cdots, x_{n}, x) \in \mathcal{O}_{n+1}(x)}\exp{\big(\sum_{i=1}^{n}\phi(x_{i},x_{i+1})\big)},$$where $\phi$ is the function defined in $(\ref{volume-m})$.

Note that
for any $\underline{x}=(x_{1},\cdots, x_{n}, x) \in \mathcal{O}_{n+1}(x)$, there exist $f_{j_{1}},\cdots, f_{j_{n}} \in \{f_{1},\cdots,f_{k}\}$ such that
$x_{i+1}=f_{j_{i}}(x_{i})$, $1 \leq {j_{i}} \leq k$, $1\leq i \leq n$. Then for any $n \in \mathbb{N}$, $x \in M$,
\begin{eqnarray*}\label{simpli}
 \Phi_{n}(x)=\sum_{\underline{x}=(x_{1},\cdots, x_{n}, x) \in \mathcal{O}_{n+1}(x)}\frac{1}{\prod_{i=1}^{n}\text{Jac}(j_{i}(x_{i}))}.
\end{eqnarray*}

By (\ref{topological pressure}), $P(T,\phi)=P(\sigma_{T},\widetilde{\phi})$.
Thus, we can calculate the equilibrium state of correspondences using the left shift map on the orbit space. However, even so, calculating the equilibrium state for the corresponding natural extension remains challenging.
When we consider the equilibrium state for the correspondence
$T$ on $M$ with finite $C^{2}$-expanding endomorphisms and the special potential $\phi$
defined in (\ref{volume-m}), we will use a type of Ruelle operator $\widetilde{\mathcal{L}}_{\phi}$, as defined in $(\ref{correspondence-operator})$, to obtain the density function  for the correspondence.
The density function is the eigenfunction of the Ruelle operator with eigenvalue $k$ (the number of the generators of $T$).
Hence we derive the  smooth probability measure on $M$ and show that it is a special kernel $\mathcal{Q}$-invariant probability measure.
We then relate the measure-theoretic entropy of the special transition probability kernel for the smooth probability measure to the skew product transformation, as described in (\ref{gene-skew measure entropy}). By applying the Abramov-Rokhlin formula (see \cite{Abramov-Rokhlin1966}) and the fiber-entropy formula  (see \cite{P.-D. Liu1999}), we obtain the equilibrium state for $P(T,\phi)$.

\begin{thmm}\label{thm:kernel-dense-function}
Let $T$ be a finitely-generated correspondence with finite $C^{2}$-expanding endomorphisms $f_{1},f_{2}$, $\cdots, f_{k}$ on $M$. Assume $\text{Coinc}(f_{j},f_{j^{\prime}})=\emptyset$ for $1 \leq j\neq j^{\prime} \leq k$. Let $\phi\colon \mathcal{O}_{2}{(T)}\to \mathbb{R}$ be the potential as in (\ref{volume-m}).
Then the equilibrium state $(\mu, \mathcal{Q})$  exists for $P(T,\phi)$. Precisely,
 \begin{align*}
  \mu=\Phi\mathbf{m} \text{  and  } \mathcal{Q}_{x}=\frac{1}{k}\sum_{j=1}^{k}\delta_{f_{j}(x)},\;\;\; \forall x \in M,
 \end{align*}where
 \begin{align*}
  \Phi(x)=\lim_{n \rightarrow +\infty}\frac{1}{k^{n}}\Phi_{n}(x),\;\;\; \forall x \in M.
 \end{align*}
\end{thmm}

\begin{remark}
  In fact, as we will see in Remark \ref{unique-sense},  the equilibrium state $(\Phi\mathbf{m},\frac{1}{k}\sum_{j=1}^{k}\delta_{f_{j}})$ for $P(T,\phi)$, obtained in Theorem \ref{thm:kernel-dense-function},  is  unique in an appropriate sense.
\end{remark}

\subsection{Bounded measurable potentials}

 Denote by $B(X,\mathbb{R})$ the set of real-valued bounded measurable
 functions on $X$. We give the following  type of topological pressure of correspondence on $X$ with respect to $\phi \in B(X,\mathbb{R})$.
\begin{definition}
Let $T$ be a forward expansive  correspondence on $X$, $\phi \in B(\mathcal{O}_{2}(T),\mathbb{R})$. Then the \emph{variational topological pressure} of $T$ with respect to potential $\phi$ is define as
\begin{equation}\label{variational topological pressure}
   P^{v}(T, \phi):= \sup_{(\mu,\mathcal{Q})}\bigg\{h_{\mu}(\mathcal{\mathcal{Q}})+ \int_{X}\int_{T(x_{1})} \phi(x_{1},x_{2})d\mathcal{Q}_{x_{1}}(x_{2})d\mu(x_{1})\bigg\},
\end{equation}
where $\mathcal{Q}$ ranges over all transition probability kernels on $X$ supported by $T$ and $\mu$
ranges over all $\mathcal{Q}$-invariant Borel probability measures on $X$. We say that $(\mu,\mathcal{Q})$ is an \emph{equilibrium state for $P^{v}(T,\phi)$} if
\begin{align*}
 P^{v}_{(\mu,\mathcal{Q})}(T, \phi)= P^{v}(T, \phi),
\end{align*}
where $P^{v}_{(\mu, \mathcal{Q})}(T, \phi)=h_{\mu}(\mathcal{\mathcal{Q}})+ \int_{X}\int_{T(x_{1})} \phi(x_{1},x_{2})d\mathcal{Q}_{x_{1}}(x_{2})d\mu(x_{1})$.
\end{definition}

If $\phi \colon \mathcal{O}_{2}(T) \to \mathbb{R}$ is continuous,
the notion of variational topological pressure $P^{v}(T, \phi)$ coincides with the  notion of topological pressure $ P(T, \phi) $ as defined in $(\ref{toppre})$. Further discussion about $P^{v}(T, \phi)$ with respect to the natural extension can be found  in Remark \ref{two-variation-coin}.

The following results show that for  the finitely-generated  correspondence of  torus systems, the variational topological pressure with respect to a real-valued bounded measurable potential is equal to $\log k$, where $k$ is the number of generators, and
 the corresponding   equilibrium states
for variational topological pressure
 exist.
\begin{thmm}\label{thmm-bounded-poten-equi}
 Let $T\colon \mathbb{T}^{m}\to \mathbb{T}^{m}$ be a  finitely-generated correspondence on the torus $\mathbb{T}^{m}$ with the generators $f_{j}=A_{j}$, $1 \leq j \leq k$, where $\{A_{j}\}_{j=1}^{k}$ are  non-singular integer matrices, and all eigenvalues of each $A_{j}$ have modulus greater than 1. Let $\phi \in B(\mathcal{O}_{2}(T),\mathbb{R})$ be defined by
\begin{align}\label{spe-poten}
 \phi(x_{1},x_{2})=\left \{\begin{array}{lcl}
                                                            -\min_{1\leq j\leq k}\{\log |\text{det}A_{j}|\}  & \mbox{ if }  & x_{1}\in E\\
                                                            -\log |\text{det}A_{j}| &\mbox{ if } & x_{1}\notin E,\text{ and } x_{2}=f_{j}(x_{1}),  1\leq j \leq k,
                                                          \end{array} \right.
\end{align}
where $E=\bigcup_{1 \leq j \neq j^{\prime} \leq k}\text{Coinc}(f_{j},f_{j^{\prime}})$.
Then the variational topological pressure $P^{v}(T,\phi)$ is equal to $\log k$ and $(\text{Leb},\sum_{j=1}^{k}\frac{1}{k}\delta_{f_{j}(x)})$, $\forall x \in \mathbb{T}^{m}$,
 is  one of the  equilibrium states for $P^{v}(T,\phi)$.
\end{thmm}

\begin{remark}
 In fact, the function $\phi \in B(\mathcal{O}_{2}(T),\mathbb{R})$ as defined in $(\ref{spe-poten})$ can be replaced by
 \begin{align*}
 \phi(x_{1},x_{2})=\left \{\begin{array}{lcl}
                                                         -c  & \mbox{ if }  & x_{1}\in E\\
                                                           -\log |\text{det}A_{j}| &\mbox{ if } & x_{1}\notin E,\text{ and } x_{2}=f_{j}(x_{1}),  1\leq j \leq k,
                                                         \end{array} \right.
\end{align*}
where $E=\bigcup_{1 \leq j \neq j^{\prime} \leq k}\text{Coinc}(f_{j},f_{j^{\prime}})$ and $c \in \{\log |\text{det}A_{1}|, \log |\text{det}A_{2}|,\cdots, \log |\text{det}A_{k}|\}$.
\end{remark}

\section{The proof of Theorem  A}

We now state some definitions from \cite{Xiaoran Li-Zhiqiang Li-Yiwei Zhang2023}.

\begin{definition}\label{strong transitivity}
  We say a correspondence $T$ on a compact connected metric space $(X, d)$ is \emph{strongly transitive} if $\bigcup_{n=1}^{+\infty}T^{-n}(x)$ is dense in $X$ for every $x \in X$.
\end{definition}

\begin{definition}\label{distance-expanding}
  We say a correspondence $T$ on a compact connected metric space $(X, d)$ is \emph{distance-expanding} if there exists $\tilde{\lambda} >1$, $\eta >0$, and $n \in \mathbb{N}$ with the property that for each $x,y \in X$ if $d(x,y) \leq \eta$, then
  \begin{equation*}\label{diatance-inf}
    \inf \{d(x^{\prime}, y^{\prime}):x^{\prime} \in T^{n}(x),\; y^{\prime} \in T^{n}(y) \}\geq \tilde{\lambda} d(x,y).
  \end{equation*}
\end{definition}

\begin{remark}
 Let $T\colon X \to X$ be a map defined by $T(x)=\{f(x)\}$, where $f\colon X\to X$ is a  single-valued continuous map. Then the definition of distance-expanding of $T$, as defined in Definition \ref{distance-expanding},  coincides with the definition of distance-expanding for single-valued continuous map $f$ (see \cite{Walters1982}, for example).
\end{remark}

Then we have the following lemmas.

\begin{lemma}\label{open-transitive-expand}
  Let $T$ be a finitely-generated correspondence with finite $C^{2}$-expanding endomorphisms $f_{1},f_{2}$, $\cdots, f_{k}$ on $M$. Assume $\text{Coinc}(f_{j},f_{j^{\prime}})=\emptyset$ for $1 \leq j\neq j^{\prime} \leq k$. Then the correspondence $T$ is strongly transitive and distance-expanding.
\end{lemma}

\begin{pf}

Since each $f_{j}$ is $C^{2}$-expanding, then for any $\varepsilon >0$, there exists $N_{j}(\varepsilon)\in \mathbb{N}$ such that
$f_{j}^{N_{j}(\varepsilon)}B(y, \varepsilon)=M, \;\forall y \in M $ (see \cite{Krzyzewski-Szlenk1969}). It follows that  $\bigcup_{n=1}^{+\infty}T^{-n}(x) \cap B(y,\varepsilon) \neq \emptyset, \;\forall x \in M$. Hence $\bigcup_{n=1}^{+\infty}T^{-n}(x)$ is dense in $M$ for every $x \in M$.
Then,
by Definition \ref{strong transitivity}, the correspondence $T$ is strongly transitive. So, in the following, we only need to prove $T$ is distance-expanding.

  Since each $f_{j}$, $1 \leq j \leq k$, is $C^{2}$-expanding, there exist positive integer $\tilde{s}_{j}$ and numbers $\varepsilon_{j}>0$, $\delta_{j} >0$, $u_{j}>1$ such that for any $x \in M$, $\#f_{j}^{-1}(x)=\tilde{s}_{j}$ (where $\tilde{s}_j$ is the degree of $f_j$), and $\min_{z\neq z^{\prime}\in f_{j}^{-1}(x)}d(z,z^{\prime}) \geq \delta_{j}$, and the preimage of the $B(x, \varepsilon_{j})$ under $f_{j}$ is a disjoint union of connected open sets $V_{i}(x,f_{j},\varepsilon_{j})$, $1 \leq i \leq \tilde{s}_{j}$, of diameter less than $\frac{2\varepsilon_{j}}{u_{j}}$, on which the restriction
\begin{equation}\label{local homeo}
  f_{j}|_{V_{i}(x,f_{j},\varepsilon_{j})}\colon V_{i}(x,f_{j},\varepsilon_{j})\to B(x,\varepsilon_{j})
\end{equation}is a  homeomorphism, and moreover, for any $y, z \in V_{i}(x,f_{j},\varepsilon_{j})$,
\begin{eqnarray}\label{one-expanding}
  d(f_{j}(y), f_{j}(z))\geq u_{j}d(y,z).
\end{eqnarray}
We now choose constants $0 < \varepsilon_{0} \leq \min_{1 \leq j \leq k}\varepsilon_{j}$, $0 < \delta_{0} \leq \min_{1 \leq j \leq k} \delta_{j}$, and $1< u_{j} \leq \min_{1\leq j \leq k}u_{j}$ such that
$$
\inf_{x \in M}\min_{1 \leq j \leq k}\min_{1 \leq i_{1}\neq i_{2} \leq \tilde{s}_{j}}d(V_{i_{1}}(x,f_{j},\varepsilon_{0}), V_{i_{2}}(x,f_{j},\varepsilon_{0}))> \delta_{0},
$$and for any $0 < \varepsilon <\varepsilon_{0}$,
\begin{equation}\label{mindiam}
\sup_{x \in M} \max_{1 \leq j \leq k} \max_{1 \leq i \leq \tilde{s}_{j}}\text{diam}V_{i}(x,f_{j},\varepsilon) \leq \frac{2\varepsilon}{u_{0}}.
\end{equation}Since $\text{Coinc}(f_{j},f_{j^{\prime}})=\emptyset$ for each pair of generators $f_{j}$ and $f_{j^{\prime}}$, $1 \leq j\neq j^{\prime} \leq k$, let
\begin{equation}\label{diandance-separeted-1}
  \varepsilon^{\star}=\min_{1\leq j \neq j^{\prime} \leq k}\inf_{x \in M}d(f_{j}(x),f_{j^{\prime}}(x)).
\end{equation}Since $M$ is compact, $\varepsilon^{\star}>0$. Then we  take $0 \leq \varepsilon^{\triangle}<\frac{\min\{\frac{\varepsilon^{\star}}{2}, \varepsilon_{0},\delta_{0}\}}{u_{0}}$, such that
\begin{equation*}\label{distance-separsted-2}
  \min_{1 \leq j \neq j^{\prime}\leq k}\inf_{x,y \in M, d(x,y)<\varepsilon^{\triangle}}d(f_{j}(x),f_{j^{\prime}}(y)) \geq \frac{\varepsilon^{\star}}{2}.
\end{equation*}
It follows that for any $1 \leq j \neq j^{\prime} \leq k$, if $d(x,y) < \varepsilon^{\triangle}< \frac{\varepsilon_{0}}{u_{0}}$, then
$$d(f_{j}(x),f_{j^{\prime}}(y))\geq \frac{\varepsilon^{\star}}{2} > \varepsilon^{\triangle}u_{0}>u_{0}d(x,y).$$
By $(\ref{local homeo})$, $(\ref{one-expanding})$ and $(\ref{mindiam})$, we also know that for any $1 \leq j \leq k$, if $d(x,y) < \varepsilon^{\triangle}$, then $d(f_{j}(x),f_{j}(y)) \geq u_{0}d(x,y)$.
So, take $\tilde{\lambda}=u_{0}>1$, $\eta=\varepsilon^{\triangle}$, and $n=1$, then for any $x,y \in M$, if $d(x,y) \leq \eta=\varepsilon^{\triangle}$, we have
\begin{align}\label{dis-expan}
  \inf\{d(x^{\prime},y^{\prime}):x^{\prime}\in T(x),\;y^{\prime}\in T(y)\} \geq \tilde{\lambda} d(x,y).
\end{align}
By definition \ref{distance-expanding}, $T$ is a distance-expanding correspondence.

This completes the proof of Lemma \ref{open-transitive-expand}.
\end{pf}

For  strongly transitive and distance-expanding correspondences, the locally H\"{o}lder continuous potentials are crucial when calculating topological pressures.
Consider a compact metric space $(X,d)$. Recall that a function $\phi\colon X\to \mathbb{R}$ is said to be  \emph{locally H\"{o}lder continuous} with an exponent $0 < \alpha < 1$ if there exist  $\xi >0$ and $C>0$ such that
\begin{align}\label{locally-continuous}
  |\phi(y)-\phi(x)|\leq C{d(x,y)}^{\alpha}
\end{align}
for all $x, y \in X$ with $d(x,y)\leq \xi$. Here, $\alpha$ is the \emph{locally H\"{o}lder exponent} and $C$ the \emph{locally H\"{o}lder constant} of $\phi$.  And $\phi$ is \emph{Hölder continuous} if the inequality  $(\ref{locally-continuous})$ holds for all $x,y \in X$.
We may assume $\xi \le \eta$, where $\eta$ as defined in Lemma \ref{open-transitive-expand}.

\begin{lemma}\label{holder potential}
  Let $T$ be a finitely-generated correspondence with finite $C^{2}$-expanding endomorphisms $f_{1},f_{2}$, $\cdots, f_{k}$ on $M$. Assume $\text{Coinc}(f_{j},f_{j^{\prime}})=\emptyset$ for $1 \leq j\neq j^{\prime} \leq k$.
  Let $\phi\colon \mathcal{O}_{2}{(T)}\to \mathbb{R}$ be the potential as in (\ref{volume-m}).
 Then $\phi$ is a locally H\"{o}lder continuous function with respect to $(\mathcal{O}_{2}{(T)}, d_{2})$.
\end{lemma}

\begin{pf}
  Since each $f_{j}$, $1 \leq j \leq k$, is $C^{2}$-expanding map on $M$, then the map $x\longmapsto -\log \text{Jac}(f_{j}(x))$ is $C^{1}$.
  Since $M$ is compact, there exists $0 <\xi \leq \eta$ such that for any  $x, y \in M$ with $d(x,y) \le \xi$, $1\leq j \leq k$,
 \begin{equation}\label{ecah-holder}
  \left|\log \text{Jac}(f_{j}(x))-\log \text{Jac}(f_{j}(y))\right|\leq C_{j}d(x,y) \leq C_{j}{d(x,y)}^{\alpha_{j}},
  \end{equation}
  where $C_{j}>0$, $0 < \alpha_{j} <1$. Let $C^{\prime}=\max_{1 \leq j \leq k}C_{j}$ and $\alpha=\min_{1 \leq j \leq k}\alpha_{j}$, $0 < \alpha <1$.
  Let $(x_1, x_2) = (x_1, f_j(x_1))$, $(y_1, y_2) = (y_1, f_{j'}(y_1))$ in $\mathcal{O}_{2}(T)$ with $d_{2}((x_1, x_2),(y_1, y_2)) \le \xi$.
  There are two cases for the indices $j$ and $j^{\prime}$.

  \emph{Case 1}. $j=j^{\prime}$. Then by $(\ref{ecah-holder})$,
 \begin{equation}\label{holder-1}
  |\phi(y_{1},y_{2})- \phi(x_{1},x_{2})|=|\log\text{Jac}(f_{j}(x_{1}))-\log\text{Jac}(f_{j}(y_{1}))| \leq C^{\prime}{d(x_{1},y_{1})}^{\alpha_{j}}.
 \end{equation}
 Since $d(x_{1},y_{1}) \le \xi$, ${d(x_{1},y_{1})}^{\alpha_{j}} \leq {d(x_{1},y_{1})}^{\alpha}$. Then $|\phi(y_{1},y_{2})- \phi(x_{1},x_{2})| \leq C^{\prime}{d(x_{1},y_{1})}^{\alpha}$.

 \emph{Case 2}. $j\neq j^{\prime}$. Then we  obtain
 \begin{align*}
  |\phi(y_{1},y_{2})- \phi(x_{1},x_{2})|=&|\log\text{Jac}(f_{j}(x_{1}))-\log\text{Jac}(f_{j^{\prime}}(y_{1}))| \\
  \leq&  |\log\text{Jac}(f_{j}(x_{1}))-\log\text{Jac}(f_{j}(y_{1}))| +|\log\text{Jac}(f_{j}(y_{1}))-\log\text{Jac}(f_{j^{\prime}}(y_{1}))|\\
   \leq &C^{\prime}{d(x_{1},y_{1})}^{\alpha} + |\log\text{Jac}(f_{j}(y_{1}))-\log\text{Jac}(f_{j^{\prime}}(y_{1}))|.
 \end{align*}

  Since $\text{Coinc}(f_{j},f_{j^{\prime}})=\emptyset$ for each pair of generators $f_{j}$ and $f_{j^{\prime}}$, $1 \leq j \neq j^{\prime} \leq k$, by (\ref{diandance-separeted-1}), we  know that $d(f_{j}(y_{1}), f_{j^{\prime}}(y_{1})) \geq \varepsilon^{\star}$, $j \neq j^{\prime}$.
  Let $$L=2 \max_{1 \leq j \leq k}\sup_{x \in M}\log \text{Jac}(f_{j}(x)),$$and
  \begin{align}\label{distancebecome}
    L^{\prime}=\max_{1 \leq j \leq k} \sup_{x \in M} \text{Jac}(f_{j}(x)).
  \end{align}

  Then we have
   \begin{align*}
 |\log\text{Jac}(f_{j}(y_{1}))-\log\text{Jac}(f_{j^{\prime}}(y_{1}))|
  \leq &\frac{L}{(\varepsilon^{\star})^{\alpha}}(\varepsilon^{\star})^{\alpha}
  \leq \frac{L}{(\varepsilon^{\star})^{\alpha}}{d(f_{j}(y_{1}), f_{j^{\prime}}(y_{1}))}^{\alpha}\\
  \leq& \frac{L}{(\varepsilon^{\star})^{\alpha}}{\big(d(f_{j}(y_{1}),f_{j}(x_{1}))+d(f_{j}(x_{1}), f_{j^{\prime}}(y_{1}))\big)}^{\alpha}\\
  \leq &\frac{L}{(\varepsilon^{\star})^{\alpha}} {\big( L^{\prime}d(x_{1},y_{1})+d(x_{2},y_{2})\big)}^{\alpha}\\
  \leq &\frac{L}{(\varepsilon^{\star})^{\alpha}}(L^{\prime}+1)^{\alpha}{d_{2}((x_{1},x_{2}),(y_{1},y_{2}))}^{\alpha}.
   \end{align*}It follows that
  \begin{equation}\label{Case-2}
   |\phi(y_{1},y_{2})- \phi(x_{1},x_{2})| \leq  \big( C^{\prime}+\frac{L}{(\varepsilon^{\star})^{\alpha}}(L^{\prime}+1)^{\alpha}\big){d_{2}((x_{1},x_{2}),(y_{1},y_{2}))}^{\alpha}.
  \end{equation}
  Let $C= C^{\prime}+\frac{L}{(\varepsilon^{\star})^{\alpha}}(L^{\prime}+1)^{\alpha}>0$. Combining (\ref{holder-1}) and $(\ref{Case-2})$,
  we obtain that for any $(x_{1},x_{2}), \; (y_{1},y_{2}) \in \mathcal{O}_{2}(T)$, if $d_{2}((x_1, x_2),(y_1, y_2)) \le \xi$, then
   $$
   |\phi(y_{1},y_{2})- \phi(x_{1},x_{2})| \leq C{d_{2}((x_{1},x_{2}),(y_{1},y_{2}))}^{\alpha},
   $$where $\alpha= \min_{1 \leq j \leq k}\alpha_{j}, \;0 < \alpha <1$.
  Hence,  $\phi\colon \mathcal{O}_{2}{(T)}\to \mathbb{R}$ is a locally  $\alpha$-H\"{o}lder continuous function with respect to $(\mathcal{O}_{2}{(T)}, d_{2})$.

  This completes the proof of Lemma \ref{holder potential}.
\end{pf}

\begin{remark}\label{extend-holder}
  From the definitions of $d_{2}$ and $d_{\omega}$, we know that if $\phi\colon \mathcal{O}_{2}{(T)}\to \mathbb{R}$ is  a locally  $\alpha$-H\"{o}lder continuous function with respect to the metric $d_{2}$ on $\mathcal{O}_{2}{(T)}$, then the function $\widetilde{\phi}\colon \mathcal{O}_{\omega}{(T)}\to \mathbb{R}$  defined in $(\ref{local constant})$ is also locally $\alpha$-H\"{o}lder continuous with respect to the metric $d_{\omega}$ on $\mathcal{O}_{\omega}{(T)}$.
\end{remark}

Using the knowledge of Riemannian geometry, the following lemma is evident.

\begin{lemma}\label{degree}
  Let $f$ be a $C^{2}$-expanding endomorphism on $(M,\mathbf{m})$. Then we have
  \begin{align*}
    \int_{M} \sum_{x \in f^{-1}(y)} \frac{1}{\text{Jac}(f(x))}d\mathbf{m}(y)=\mathbf{m}(M).
  \end{align*}
\end{lemma}







Let $T$ be a finitely-generated correspondence with finite $C^{2}$-expanding endomorphisms $f_{1},f_{2}, \cdots, f_{k}$ on $M$. Assume $\text{Coinc}(f_{j},f_{j^{\prime}})=\emptyset$ for $1 \leq j\neq j^{\prime} \leq k$. Let $\phi\colon \mathcal{O}_{2}{(T)}\to \mathbb{R}$ be a function defined as in $(\ref{volume-m})$.
In the following, we will consider some backward orbits under the  correspondence $T$.
Recall that
\begin{align*}
 \mathcal{O}_{n+1}(x)=&\{\underline{x}=(x_{1},\cdots, x_{n}, x_{n+1}) \in \mathcal{O}_{n+1}(T):
 x_{n+1}=x\},
\end{align*}and
$$\Phi_{n}(x)=\sum_{\underline{x}=(x_{1},\cdots, x_{n}, x) \in \mathcal{O}_{n+1}(x)}\exp{\big(\sum_{i=1}^{n}\phi(x_{i},x_{i+1})\big)}.$$
Let \begin{align*}
 \mathcal{O}_{n+1}(x,(j_{1},\cdots,j_{n})):=\{\underline{x}=(x_{1},\cdots, x_{n}, x) \in \mathcal{O}_{n+1}(x):
  x_{i+1}=f_{j_{i}}(x_{i}), 1\leq i \leq n\}.
\end{align*}

 Then we have the following Lemma \ref{bounded-control} and Lemma \ref{log-calculate}.

\begin{lemma}\label{bounded-control}
   There exists a constant $\widetilde{C}>0$, such that for any $x, y \in M$, $n \in \mathbb{N}$,
   \begin{equation*}
     \frac{1}{\widetilde{C}} \leq \frac{\Phi_{n}(x)}{\Phi_{n}(y)} \leq \widetilde{C}.
   \end{equation*}
\end{lemma}

\begin{pf}
  By Lemma \ref{holder potential},  we know that the function $\phi \colon \mathcal{O}_{2}{(T)}\to \mathbb{R}$ defined in $(\ref{volume-m})$ is locally H\"{o}lder with respect to the metric $d_{2}$ on $\mathcal{O}_{2}{(T)}$.
  Define $\xi$, $\alpha$ and $C$ as in Lemma \ref{holder potential}.  Note that
if $d(x, y) \leq  \xi$, then for any element $\underline{x}=(x_{1},\cdots, x_{n}, x) \in \mathcal{O}_{n+1}(x,(j_{1},\cdots,j_{n}))$, there is a unique element $\underline{y}=(y_{1},\cdots, y_{n}, y) \in
\mathcal{O}_{n+1}(y,(j_{1},\cdots,j_{n}))$ such that $d_{n+1}(\underline{x},\underline{y})=d(x,y)\leq \xi$, and vice versa.
By $(\ref{dis-expan})$, for such a pair of elements $\underline{x}$ and $\underline{y}$, we have $d(x_{i+1},y_{i+1}) \leq {\tilde{\lambda}}^{-(n-i)}\xi$ for every $1 \leq i \leq n$. Hence,
  $$
 \big|\sum_{i=1}^{n}\phi(x_{i},x_{i+1})-\sum_{i=1}^{n}\phi(y_{i},y_{i+1})\big| \leq \frac{C}{1-{\tilde{\lambda}}^{-\alpha}}\xi^{\alpha}.
  $$
  Thus, we have
  $$
  \frac{1}{\exp{\big(\frac{C}{1-{\tilde{\lambda}}^{-\alpha}}\xi^{\alpha}\big)}} \leq \frac{\Phi_{n}(x)}{\Phi_{n}(y)}=\frac{\sum_{j_{1},\cdots,j_{n}=1}^{k} \sum_{\underline{x}\in \mathcal{O}_{n+1}(x,(j_{1},\cdots,j_{n}))}\exp{\big(\sum_{i=1}^{n}\phi(x_{i},x_{i+1})\big)}}{\sum_{j_{1},\cdots,j_{n}=1}^{k} \sum_{\underline{y}\in \mathcal{O}_{n+1}(y,(j_{1},\cdots,j_{n}))}\exp{\big(\sum_{i=1}^{n}\phi(y_{i},y_{i+1})\big)}} \leq \exp{\big(\frac{C}{1-{\tilde{\lambda}}^{-\alpha}}\xi^{\alpha}\big)}.
  $$

  For $\xi >0$, there exists $N_{j}(\xi)\in \mathbb{N}$
 such that
$f_{j}^{N_{j}(\xi)}B(x, \xi)=M$, $\forall x \in M$. Take $N(\xi):=\max_{1 \leq j \leq k}N_{j}(\xi)$. Then for any $x,y \in M$, there exists $0 \leq m(\xi)\leq  N(\xi)$ such that $(T^{m(\xi)}B(x,\xi))\cap B(y,\xi)\neq  \emptyset$.

  Now fix arbitrary $x,y \in M$. Then there exists a point $y^{\prime} \in  B(x, \xi)$ such that $y^{\prime\prime} \in T^{m(\xi)}(y^{\prime})$ and $y^{\prime\prime} \in B(y, \xi)$. Thus, we have
  \begin{align*}
    \Phi_{n}(x)=& \sum_{\underline{x}=(x_{1},\cdots, x_{n}, x) \in \mathcal{O}_{n+1}(x)}\exp{\big(\sum_{i=1}^{n}\phi(x_{i},x_{i+1})\big)}\\
    \leq & \exp{\big(\frac{C}{1-{\tilde{\lambda}}^{-\alpha}}\xi^{\alpha}\big)}\big(\sum_{\underline{y^{\prime}}=(y^{\prime}_{1},\cdots, y^{\prime}_{n}, y^{\prime}) \in \mathcal{O}_{n+1}(y^{\prime})}\exp{\big(\sum_{i=1}^{n}\phi(y^{\prime}_{i},y^{\prime}_{i+1})\big)} \big)\\
    \leq & \exp{\big(\frac{C}{1-{\tilde{\lambda}}^{-\alpha}}\xi^{\alpha}\big)}(L^{\prime})^{m(\xi)} \big( \sum_{\underline{y^{\prime\prime}}=(y^{\prime\prime}_{1},\cdots, y^{\prime\prime}_{m(\xi)+n},y^{\prime\prime})\in \mathcal{O}_{m(\xi)+n+1}(y^{\prime\prime})}  \exp{\big(\sum_{i=1}^{m(\xi)+n}\phi(y^{\prime\prime}_{i},y^{\prime\prime}_{i+1})\big)} \big)\\
    \leq & \exp{\big(\frac{C}{1-{\tilde{\lambda}}^{-\alpha}}\xi^{\alpha}\big)}(L^{\prime})^{m(\xi)}\big(\sum_{j=1}^{k}\text{deg}(f_{j})\big)^{m(\xi)}   \big( \sum_{\underline{y^{\prime\prime}}=(y^{\prime\prime}_{1},\cdots, y^{\prime\prime}_{n},y^{\prime\prime})\in \mathcal{O}_{n+1}(y^{\prime\prime})}  \exp{\big(\sum_{i=1}^{n}\phi(y^{\prime\prime}_{i},y^{\prime\prime}_{i+1})\big)} \big)  \\
    \leq & \exp{\big(\frac{2C}{1-{\tilde{\lambda}}^{-\alpha}}\xi^{\alpha}\big)}(L^{\prime})^{m(\xi)}\big(\sum_{j=1}^{k}\text{deg}(f_{j})\big)^{m(\xi)}   \big( \sum_{\underline{y}=(y_{1},\cdots, y_{n},y)\in \mathcal{O}_{n+1}(y)}  \exp{\big(\sum_{i=1}^{n}\phi(y_{i},y_{i+1})\big)} \big) \\
    = &\exp{\big(\frac{2C}{1-{\tilde{\lambda}}^{-\alpha}}\xi^{\alpha}\big)}(L^{\prime})^{m(\xi)}\big(\sum_{j=1}^{k}\text{deg}(f_{j})\big)^{m(\xi)} \Phi_{n}(y),
  \end{align*}where $L^{\prime}$ is defined as in $(\ref{distancebecome})$.

  It follows that there exists a constant $\widetilde{C}:=\exp{\big(\frac{2C}{1-{\tilde{\lambda}}^{-\alpha}}\xi^{\alpha}\big)}(L^{\prime})^{m(\xi)}(\sum_{j=1}^{k}\text{deg}(f_{j}))^{m(\xi)}>0$ such that for any $x, y \in M$, $n \in \mathbb{N}$,
  $$
  \frac{1}{\widetilde{C}} \leq \frac{\Phi_{n}(x)}{\Phi_{n}(y)} \leq \widetilde{C}.
  $$

  This completes the proof of Lemma \ref{bounded-control}.
\end{pf}

By Lemma \ref{holder potential}, Lemma \ref{degree} and Lemma \ref{bounded-control}, we have the following results.

\begin{lemma}\label{log-calculate}
  For any $n \in \mathbb{N}$, $x \in M$,
  \begin{equation*}
   \lim_{n\rightarrow +\infty}\frac{1}{n}\log \Phi_{n}(x)=\log k.
  \end{equation*}
\end{lemma}

\begin{pf}

 By Lemma \ref{degree}, we have
  \begin{align}\label{apply-mvol}
 \int_{M} \Phi_{n}(x) d\mathbf{m}(x)
  = &  \sum_{j_{1},\cdots,j_{n}=1}^{k}  \int_{M}  \sum_{\underline{x}\in \mathcal{O}_{n+1}(x,(j_{1},\cdots,j_{n}))} \frac{1}{\prod_{i=1}^{n}   \text{Jac}(f_{j_{i}}(x_{i}))} d\mathbf{m}(x) \notag
   \\ =& \sum_{j_{1},\cdots,j_{n}=1}^{k}  \int_{M} \sum_{x_{1} \in (f_{j_{n}}\circ\cdots\circ f_{j_{1}})^{-1}(x)}\frac{1}{\text{Jac}((f_{j_{n}}\circ\cdots\circ f_{j_{1}})(x_{1}))}   d\mathbf{m}(x) \notag
   \\=&\mathbf{m}(M)\cdot k^{n}.
  \end{align}

  We claim:
   For any $n \in \mathbb{N}$, $\Phi_{n}\colon M\to \mathbb{R}$ is a  continuous function on $M$.
In fact, we will prove it is  H\"{o}lder continuous.

  For any $x, y \in M$, we have
    \begin{align*}
      &|\Phi_{n}(x)-\Phi_{n}(y)|  \leq  \\&
      \sum_{j_{1},\cdots,j_{n}=1}^{k} \bigg|\sum_{\underline{x}\in \mathcal{O}_{n+1}(x,(j_{1},\cdots,j_{n}))} \frac{1}{\prod_{i=1}^{n}\text{Jac}(f_{j_{i}}(x_{i}))}-\sum_{\underline{y}\in \mathcal{O}_{n+1}(y,(j_{1},\cdots,j_{n}))} \frac{1}{\prod_{i=1}^{n}\text{Jac}(f_{j_{i}}(y_{i}))}\bigg|.
    \end{align*}
    Define $\eta$ as in Lemma \ref{open-transitive-expand}. Similar to what we have done in the proof of Lemma \ref{bounded-control},
if $d(x, y) \leq \eta$, then for any element $\underline{x}=(x_{1},\cdots, x_{n}, x) \in \mathcal{O}_{n+1}(x,(j_{1},\cdots,j_{n}))$, there is a unique element $\underline{y}=(y_{1},\cdots, y_{n}, y) \in
\mathcal{O}_{n+1}(y,(j_{1},\cdots,j_{n}))$ such that $d_{n+1}(\underline{x},\underline{y})=d(x,y)\leq \eta$, and vice versa.
Thus, for such a pair of elements $\underline{x}$ and $\underline{y}$, we have
\begin{align*}
  \bigg|  \frac{1}{\prod_{i=1}^{n}\text{Jac}(f_{j_{i}}(x_{i}))}- \frac{1}{\prod_{i=1}^{n}\text{Jac}(f_{j_{i}}(y_{i}))}\bigg|
  \leq & \bigg|   \frac{\prod_{i=1}^{n}\text{Jac}(f_{j_{i}}(x_{i}))-\prod_{i=1}^{n}\text{Jac}(f_{j_{i}}(y_{i}))}{\tilde{L}^{2n}} \bigg| \\
  \leq & \frac{nL^{\prime}Cd(x,y)^{\alpha}}{\tilde{L}^{2n}},
\end{align*}where $\tilde{L}=\min_{1\leq j \leq k} \inf_{x \in M} \text{Jac}(f_{j}(x))$, $C$, $L^{\prime}$, and $\alpha$ are as defined  in Lemma \ref{holder potential}.

Note that $\#\mathcal{O}_{n+1}(x,(j_{1},\cdots,j_{n}))=\prod_{i=1}^{n}\text{deg}(f_{j_{i}})$. Then for any $x,y \in M$, if $d(x, y) \leq \eta$, we have
\begin{align*}
 &|\Phi_{n}(x)-\Phi_{n}(y)| \leq \sum_{j_{1},\cdots,j_{n}=1}^{k}\prod_{i=1}^{n}\text{deg}(f_{j_{i}}) \frac{nL^{\prime}C}{\tilde{L}^{2n}}d(x,y)^{\alpha}=\frac{nL^{\prime}C\big(\sum_{j=1}^{k} \text{deg}(f_{j})\big)^{n}}{\tilde{L}^{2n}}d(x,y)^{\alpha}.
\end{align*}
And for any $x, y \in M$, if $d(x,y)> \eta$, then
$$
|\Phi_{n}(x)-\Phi_{n}(y)| \leq
\frac{2{L^{\prime}}^{n}\big(\sum_{j=1}^{k} \text{deg}(f_{j})\big)^{n}}{\tilde{L}^{2n}{\eta}^{\alpha}}{\eta}^{\alpha}
\leq\frac{2{L^{\prime}}^{n}\big(\sum_{j=1}^{k} \text{deg}(f_{j})\big)^{n}}{\tilde{L}^{2n}{\eta}^{\alpha}}d(x,y)^{\alpha}.
$$
Take $\widetilde{C_{n}}=\max\Big\{\frac{nL^{\prime}C\big(\sum_{j=1}^{k} \text{deg}(f_{j})\big)^{n}}{\tilde{L}^{2n}},\frac{2{L^{\prime}}^{n}\big(\sum_{j=1}^{k} \text{deg}(f_{j})\big)^{n}}{\tilde{L}^{2n}{\eta}^{\alpha}}\Big\}$. Thus, $\Phi_{n}\colon M\to \mathbb{R}$ is  an $\alpha$-H\"{o}lder continuous function on $M$ with H\"{o}lder exponent $\alpha$ ($0 < \alpha <1$) and H\"{o}lder constant $\widetilde{C_{n}}$.
This completes the proof of the Claim.

By Lemma \ref{bounded-control} and the H\"{o}lder continuity property of $\Phi_{n}\colon M\to \mathbb{R}$, we have
\begin{equation*}\label{max-C-min}
  \max_{x \in M}\Phi_{n}(x) \leq \widetilde{C}\min_{x\in M}\Phi_{n}(x).
\end{equation*}
And by $(\ref{apply-mvol})$, we have
\begin{equation*}\label{min-mn}
  (\min_{x\in M}\Phi_{n}(x))\cdot\mathbf{m}(M)=\int_{M}\min_{x\in M}\Phi_{n}(x) d\mathbf{m}(x) \leq \int_{M} \Phi_{n}(x) d\mathbf{m}(x) =\mathbf{m}(M)\cdot k^{n}.
\end{equation*}
Thus, we have $\min_{x\in M}\Phi_{n}(x) \leq k^{n}$ and $\max_{x \in M}\Phi_{n}(x) \leq \widetilde{C}k^{n}$. Similarly, we also have  $\max_{x\in M}\Phi_{n}(x) \geq k^{n}$ and $\min_{x \in M}\Phi_{n}(x) \geq \frac{1}{\widetilde{C}}k^{n}$.
Consequently,
$$
\frac{1}{\widetilde{C}}k^{n} \leq \min_{x\in M}\Phi_{n}(x) \leq k^{n} \leq \max_{x\in M}\Phi_{n}(x) \leq \widetilde{C}k^{n},
$$which implies for any $x \in M$,
$$
\frac{1}{\widetilde{C}}k^{n} \leq  \Phi_{n}(x) \leq \widetilde{C}k^{n},
$$where $\widetilde{C}$ is the constant independent of any $x \in M$ and $n \in \mathbb{N}$.

It follows that
$$
\lim_{n\rightarrow +\infty} \frac{1}{n}\log \Phi_{n}(x)=\log k.
$$

This completes the proof of Lemma \ref{log-calculate}.
\end{pf}

 Now we are ready to prove Theorem \ref{thm:pressure-m}.

\begin{potA}

Let $K^{\prime}(\varepsilon)=\{x^{(l)} \in M: 1 \leq l \leq K(\varepsilon)\}$ be an $\varepsilon$-net of $M$.
Recall that $\mathcal{O}_{n+1}(x^{(l)})=\{\underline{x}=(x_{1},\cdots,x_{n},x^{(l)}) \in \mathcal{O}_{n+1}(T)\}$. We claim that the set
 $$
 \mathcal{O}_{n+1}^{\prime}(T)=\{\underline{x}=(x_{1},\cdots,x_{n},x^{(l)}) \in \mathcal{O}_{n+1}(T): x^{(l)} \in K^{\prime}(\varepsilon)\}=\bigcup_{x^{(l)} \in K^{\prime}(\varepsilon)}\mathcal{O}_{n+1}(x^{(l)})
 $$is a $(d_{n+1},\varepsilon)$-spanning set of $(\mathcal{O}_{n+1}(T),d_{n+1})$.

 In fact, for any $\underline{y}=(y_{1},y_{2},\cdots,y_{n+1}) \in \mathcal{O}_{n+1}(T)$ with $y_{i+1}=f_{j_{i}}(y_{i})$, $1 \leq j_{i} \leq k$, $1 \leq i \leq n $, we  choose $\underline{x}=(x_{1},\cdots,x_{n},x^{(l)}) \in \mathcal{O}_{n+1}^{\prime}(T)$ such that
 $d(x^{(l)}, y_{n+1}) < \varepsilon$, $x_{i+1}=f_{j_{i}}(x_{i})$, $1 \leq i \leq n-1$, $x^{(l)}=f_{j_{n}}(x_{n})$, and $d(x_{i}, y_{i}) \leq \frac{\varepsilon}{{\tilde{\lambda}}^{n+1-i}} < \varepsilon$, $1 \leq i \leq n$, $\tilde{\lambda}$ is defined as in $(\ref{dis-expan})$.
By applying  Lemma \ref{bounded-control}, we have
\begin{align}\label{topspagene}
   P(T,\phi) \leq& \lim\limits_{\varepsilon\rightarrow 0} \limsup_{n\rightarrow +\infty}\frac{1}{n}\log\big\{ \big(\max_{x^{(l)}\in K^{\prime}(\varepsilon)} \Phi_{n}(x^{(l)}) \big)\cdot K(\varepsilon)\big\}
  \leq \limsup_{n \rightarrow +\infty}\frac{1}{n}\log (\widetilde{C}\Phi_{n}(x)) \notag \\
  =& \lim_{n \rightarrow +\infty}\frac{1}{n}\log \Phi_{n}(x), \;\;\;\forall x \in M.
\end{align}

Note that if  $\text{Coinc}(f_{j}, f_{j^{\prime}})=\emptyset$ for  each pair of generators $f_{j}$ and $f_{j^{\prime}}$, $1\leq  j\neq j^{\prime}\leq k$, then take $\eta$ as in Lemma \ref{open-transitive-expand}. For any $x^{(l)} \in K^{\prime}(\varepsilon)$, the set $\mathcal{O}_{n+1}(x^{(l)})$ is a $(d_{n+1},\eta)$-separated set of $(\mathcal{O}_{n+1}(T),d_{n+1})$. Then by the definition of $P(T,\phi)$ and Lemma \ref{bounded-control}, we  have
\begin{align}\label{onebranchspan}
 P(T,\phi) \geq \lim_{n \rightarrow +\infty}\frac{1}{n} \log \Phi_{n}(x),\;\;\; \forall x \in M.
\end{align}

By $(\ref{topspagene})$ and $(\ref{onebranchspan})$, we know that for any $x \in M$,
  \begin{equation*}\label{orbit-top-pre}
    P(T,\phi)=\lim_{n\rightarrow +\infty}\frac{1}{n} \log \Phi_{n}(x).
  \end{equation*}
Thus, by Lemma \ref{log-calculate}, we have $P(T,\phi)=\log k$.

  This completes the proof of Theorem \ref{thm:pressure-m}.
\end{potA}

\section{The proof of Theorem B}
In this section, our target is to prove Theorem \ref{thm:kernel-dense-function}. The key tool for our proof is the Ruelle operator for the finitely-generated correspondence $T$ and the locally H\"{o}lder potential $\phi$.

Let $f\colon X\to X$ be a continuous map on the compact metric space $(X,d)$. Let $\psi$ be a continuous function on $X$. We recall the definition of the Ruelle operator $\mathcal{L}_{\psi}$ on continuous functions on $X$ given by
\begin{equation*}\label{Ruelle-operator}
  \mathcal{L}_{\psi}(\varphi)(x):=\sum_{y\in f^{-1}(x)}\exp(\psi(y))\cdot\varphi(y).
\end{equation*}

Ruelle operator plays an important role in finding the invariant measures for certain systems, including smooth systems. Krzyzewski and Szlenk \cite{Krzyzewski-Szlenk1969} considered the existence and the uniqueness of a smooth invariant  measure for a $C^{2}$-expanding map on $M$. Precisely, if $f\colon M\to M$ is a $C^{2}$-expanding map, then there exists an $f$-invariant probability  measure  $\mu$ that is absolutely continuous with respect to $\mathbf{m}$. If we denote $\rho^{\prime}=\frac{d\mu}{d \mathbf{m}}$, then $\mu =\rho^{\prime}\mathbf{m}$ and
\begin{align*}\label{dense-single-expanding}
 \rho^{\prime}(x)=\lim_{n \rightarrow +\infty} \sum_{y \in f^{-n}(x)} \frac{1}{\text{Jac}(f^{n}(y))}, \;\;\; \forall x \in M.
\end{align*}
Let $\psi(x)=-\log \text{Jac}(f(x))$, $\forall x \in M$.
For any positive continuous function $\varphi$ on $M$, denote
\begin{align*}\label{fixedpoint}
  \mathcal{L}_{\psi}(\varphi)(x):=\sum_{y \in f^{-1}(x)}\frac{\varphi(y)}{\text{Jac}(f(y))},
\end{align*}
then $\mathcal{L}_{\psi}(\rho^{\prime})(x)=\rho^{\prime}(x)$, $\forall x \in M$.

Let $T$ be a finitely-generated correspondence with finite $C^{2}$-expanding endomorphisms $f_{1},f_{2},\cdots, f_{k}$ on $M$. Assume $\text{Coinc}(f_{j},f_{j^{\prime}})=\emptyset$ for $1 \leq j\neq j^{\prime} \leq k$.
Inspired by  the work of Krzyzewski and Szlenk \cite{Krzyzewski-Szlenk1969},
we then  obtain the following results for the correspondence $T$ and the potential $\phi$ under the setting of Theorem \ref{thm:kernel-dense-function}.

For any positive continuous function $\varphi$ on $M$, denote
\begin{align}\label{correspondence-operator}
 \widetilde{\mathcal{L}}_{\phi}(\varphi)(x):=\sum_{(x_{1},x)\in \mathcal{Q}_{2}(x)}\exp(\phi(x_{1},x))\cdot\varphi(x_{1})=\sum_{(x_{1},x)\in \mathcal{Q}_{2}(x)}\frac{\varphi(x_{1})}{\text{Jac}(f_{j_{1}}(x_{1}))}, \text{ if }x=f_{j_{1}}(x_{1}),\;\forall x \in M.
\end{align}
Then for every $x \in M$, we have
\begin{align}\label{corresp-fixed}
  \widetilde{\mathcal{L}}_{\phi}(\Phi)(x)=k\Phi(x),
\end{align}where
\begin{align}\label{dense-single-expanding}
 \Phi(x)=\lim_{n \rightarrow +\infty}\frac{1}{k^{n}}\sum_{\underline{x}=(x_{1},\cdots, x_{n}, x) \in \mathcal{O}_{n+1}(x)}\exp{\big(\sum_{i=1}^{n}\phi(x_{i},x_{i+1})\big)}=\lim_{n \rightarrow +\infty}\frac{1}{k^{n}}\Phi_{n}(x), \;\;\;\forall x \in M.
\end{align}

\begin{remark}
  In particular, if $k=1$, then $\Phi(x)=\rho^{\prime}(x)$, and $\widetilde{\mathcal{L}}_{\phi}(\Phi)(x)=\mathcal{L}_{\psi}(\rho^{\prime})(x)$, $\forall x \in M$.
\end{remark}

Now we are ready to prove Theorem \ref{thm:kernel-dense-function}.

\begin{potB}

\emph{Step 1.}  Show that $(\Phi\mathbf{m})\mathcal{Q}=\Phi\mathbf{m}$, where $\Phi(x)$ as defined in $(\ref{dense-single-expanding})$ and $\mathcal{Q}_{x}=\frac{1}{k}\sum_{j=1}^{k}\delta_{f_{j}}(x), \forall x \in M$.

By Lemma \ref{open-transitive-expand}, we know that the correspondence $T$ is distance-expanding with the expansive constant $\eta$. Then for any $\frac{\eta}{2}$ ball $B(x,\frac{\eta}{2})$, $\forall x \in M$, for any $1\leq j \leq k$, $f_{j}$ on each component $V_{i}(x,f_{j},\frac{\eta}{2})$ of $f_{j}^{-1}(B(x,\frac{\eta}{2}))$, $1 \leq i \leq \text{deg}(f_{j})$, is injective and the restriction $f_{j}|_{V_{i}(x,f_{j},\frac{\eta}{2})}\colon V_{i}(x,f_{j},\frac{\eta}{2})\to B(x,\frac{\eta}{2})$ is a homeomorphism.
Furthermore, $V_{i}(x,f_{j},\frac{\eta}{2})\cap V_{i^{\prime}}(x,f_{j^{\prime}},\frac{\eta}{2})=\emptyset$, $1 \leq j \neq j^{\prime}\leq k$, $1 \leq i \leq \text{deg}(f_{j})$, $1 \leq i^{\prime} \leq \text{deg}(f_{j^{\prime}})$.

By $(\ref{correspondence-operator})$ and $(\ref{corresp-fixed})$, we know that $\widetilde{\mathcal{L}}_{\phi}\Phi(x):=\sum_{(x_{1},x)\in \mathcal{O}_{2}(x)}\frac{\Phi(x_{1})}{\text{Jac}(f_{j_{1}}(x_{1}))}$, $\forall x \in M$ and $\widetilde{\mathcal{L}}_{\phi}\Phi(x)=k\Phi(x)$. Hence, for any $B(x,\frac{\eta}{2})$, $\forall x \in M$,
\begin{align*}
   (\Phi\mathbf{m})(B(x,\frac{\eta}{2})) =& \int_{B(x,\frac{\eta}{2})} \Phi(y) d\mathbf{m}(y)=\int_{B(x,\frac{\eta}{2})} \frac{1}{k}\widetilde{\mathcal{L}}_{\phi}\Phi(y) d\mathbf{m}(y)  \notag
   \\=&  \frac{1}{k}\sum_{(y_{1},y)\in \mathcal{O}_{2}(y)}\int_{B(x,\frac{\eta}{2})} \frac{\Phi(y_{1})}{\text{Jac}(f_{j_{1}}(y_{1}))}      d\mathbf{m}(y) \notag
   \\=& \frac{1}{k}\sum_{1 \leq j_{1}\leq k}\sum_{1 \leq i \leq \text{deg}(f_{j_{1}})}\int_{V_{i}(x,f_{j_{1}},\frac{\eta}{2})}\Phi(y_{1})     d\mathbf{m}(y_{1})                                 \notag
  \\=& \frac{1}{k}\sum_{1 \leq j_{1}\leq k}(\Phi\mathbf{m})(f_{j_{1}}^{-1}B(x,\frac{\eta}{2})) \notag   \\
  =& ((\Phi\mathbf{m})\mathcal{Q})(B(x,\frac{\eta}{2})).
\end{align*}
Hence, we get $(\Phi\mathbf{m})\mathcal{Q}=\Phi\mathbf{m}$.

\emph{Step 2.}  Prove that $P(T,\phi)=h_{(\Phi\mathbf{m})}(\mathcal{Q})+\int_{M}\int_{T(x_{1})} \phi(x_{1}, x_{2})   d\mathcal{Q}_{x_{1}}(x_{2}) d(\Phi\mathbf{m})(x_{1})$.

Note that $\int_{M}\int_{T(x_{1})} \phi(x_{1}, x_{2})   d\mathcal{Q}_{x_{1}}(x_{2}) d(\Phi\mathbf{m})(x_{1})=\frac{1}{k}\sum_{j=1}^{k}\int_{M} \log \text{Jac}(f_{j}(x))    d(\Phi\mathbf{m})(x)$. By Theorem \ref{thm:pressure-m}, $P(T,\phi)=\log k$. Hence, we only need to prove $h_{(\Phi\mathbf{m})}(\mathcal{Q})=\log k + \frac{1}{k}\sum_{j=1}^{k}\int_{M} \log \text{Jac}(f_{j}(x))    d(\Phi\mathbf{m})(x)$.

Since $(\Phi\mathbf{m})\mathcal{Q}=\Phi\mathbf{m}$, then $(\Phi\mathbf{m})\mathcal{Q}^{\omega} \in \mathcal{M}_{\sigma_{T}}(\mathcal{O}_{\omega}(T))$. By Remark \ref{forward-infite}, we  know that
\begin{eqnarray*}
  (\Phi\mathbf{m})\mathcal{Q}^{\omega}(A_{1} \times \cdots \times A_{n}\times \{ {M}^{\omega}\})=(\Phi\mathbf{m})\mathcal{Q}^{[n-1]}(A_{1} \times \cdots \times A_{n}), \forall A_{j} \in \mathcal{B}(M),\;1 \leq j \leq n.
\end{eqnarray*}Hence, by Definition \ref{parti-entropy and entropy}, we have
\begin{align}
\begin{split}\label{measure entropy}
h_{(\Phi\mathbf{m})}(\mathcal{Q})=
  \sup_{\mathcal{A}}\lim_{n \rightarrow +\infty}\frac{-1}{n}\sum_{A_{1}\times \cdots \times A_{n} \in \mathcal{A}^{n}}(\Phi\mathbf{m})\mathcal{Q}^{[n-1]}(A_{1} \times \cdots \times A_{n}) \log (\Phi\mathbf{m})\mathcal{Q}^{[n-1]}(A_{1} \times \cdots \times A_{n}),
  \end{split}
\end{align}where $\mathcal{A}$ ranges over all finite measurable partitions of $M$ and $\mathcal{A}^{n}=\{A_{1}\times \cdots \times A_{n}: A_{j} \in \mathcal{A}, 1 \leq j \leq n\}$.

Let $\Sigma_{k}=\prod_{n \in \mathbb{N}}\{1,2,\cdots,k\}$ be the standard symbolic space with the product topology. Define a map $\widetilde{\sigma}\colon \Sigma_{k} \times M \to \Sigma_{k} \times M$ by
$$
\widetilde{\sigma}(\{j_{n}\}_{n \in \mathbb{N}}, x)=(\{j_{n+1}\}_{n \in \mathbb{N}}, f_{j_{1}}(x)).
$$This is a \textit{skew product over the shift transformation} $\sigma_{k}\colon \Sigma_{k} \to \Sigma_{k}$ by $\sigma_{k}(\{j_{n}\}_{n \in \mathbb{N}})=\{j_{n+1}\}_{n \in \mathbb{N}}$.
Define a map $\widetilde{\pi}\colon \Sigma_{k} \times M \to \mathcal{O}_{\omega}(T)$ by
$$
\widetilde{\pi}((\{j_{n}\}_{n \in \mathbb{N}},x))=(x_{1},x_{2},x_{3},\cdots),
$$where $x_{1}=x$ and $x_{n}=f_{j_{n-1}} \circ \cdots \circ f_{j_{2}} \circ f_{j_{1}}(x_{1})$ for $n \geq 2$.
Then $\sigma_T\colon \mathcal{O}_{\omega}{(T)}\to \mathcal{O}_{\omega}{(T)}$ is the factor map of $\widetilde{\sigma}\colon \Sigma_{k} \times M \to \Sigma_{k} \times M$.
By the definition of $(\Phi\mathbf{m})\mathcal{Q}^{[n-1]}$, we obtain
\begin{align}\label{calculation-n-kernel}
\begin{split}
&(\Phi\mathbf{m})\mathcal{Q}^{[n-1]}(A_{1} \times \cdots \times A_{n})   \\=&
  \sum_{j_{1},j_{2},\cdots,j_{n-1}=1}^{k}\frac{1}{k^{n-1}} (\Phi\mathbf{m})(A_{1}\cap f_{j_{1}}^{-1}A_{2}\cap(f_{j_{2}}\circ f_{j_{1}})^{-1}A_{3}\cap \cdots \cap (f_{j_{n-1}}\circ \cdots \circ f_{j_{2}} \circ f_{j_{1}})^{-1}A_{n}).
  \end{split}
\end{align}Combining  $(\ref{measure entropy})$ and $(\ref{calculation-n-kernel})$, we have
\begin{eqnarray}\label{pi-product project}
 h_{(\Phi\mathbf{m})}(\mathcal{Q})=h_{\widetilde{\pi}({\mathbb{P}_{k}}^{\mathbb{N}}\times (\Phi\mathbf{m}))}(\sigma_{T}),
\end{eqnarray}where $\mathbb{P}_{k}=(\underbrace{\frac{1}{k}, \frac{1}{k}, \cdots, \frac{1}{k}}_{k})$, and ${\mathbb{P}_{k}}^{\mathbb{N}}$ is a $\sigma_{k}$-ergodic probability measure on $\Sigma_{k}$ which gives the point $j$, $1 \leq j \leq k$, measure $\frac{1}{k}$.

Since $(\Phi\mathbf{m})\mathcal{Q}=\Phi\mathbf{m}$, then $(\Phi\mathbf{m})(A)=\frac{1}{k}\sum_{1 \leq j \leq k}(\Phi\mathbf{m})(f_{j}^{-1}A)$, $\forall A \in  \mathcal{B}(M)$, and hence
${\mathbb{P}_{k}}^{\mathbb{N}}\times (\Phi\mathbf{m})$ is $\widetilde{\sigma}$-invariant.

From the  proof of the main theorem in \cite{Zhu-Zhang-Shi2014}, we  know that $\widetilde{\pi}$ is one to one on a set of full ${\mathbb{P}_{k}}^{\mathbb{N}}\times (\Phi\mathbf{m})$ measure.
Hence
\begin{eqnarray}\label{enptyset-product-entroy}
 h_{{\mathbb{P}_{k}}^{\mathbb{N}}\times (\Phi\mathbf{m})}(\widetilde{\sigma})=h_{\widetilde{\pi}({\mathbb{P}_{k}}^{\mathbb{N}}\times (\Phi\mathbf{m}))}(\sigma_{T}).
\end{eqnarray}
 Combining $(\ref{pi-product project})$ and $(\ref{enptyset-product-entroy})$, we obtain
 \begin{eqnarray}\label{gene-skew measure entropy}
  h_{(\Phi\mathbf{m})}(\mathcal{Q})= h_{{\mathbb{P}_{k}}^{\mathbb{N}}\times (\Phi\mathbf{m})}(\widetilde{\sigma}).
 \end{eqnarray}
 By Abramov-Rokhlin's formula (see \cite{Abramov-Rokhlin1966}, for example), we have
 \begin{eqnarray}\label{gene-fiber}
   h_{{\mathbb{P}_{k}}^{\mathbb{N}}\times (\Phi\mathbf{m})}(\widetilde{\sigma})= h_{{\mathbb{P}_{k}}^{\mathbb{N}}}(\sigma_{k})+h_{(\Phi\mathbf{m})}(f),
 \end{eqnarray}where $h_{{\mathbb{P}_{k}}^{\mathbb{N}}}(\sigma_{k})=\log k$ and $h_{(\Phi\mathbf{m})}(f)$ denotes the \emph{fiber-entropy} of $\widetilde{\sigma}$ with respect to the measure $\Phi\mathbf{m}$. Then by Theorem $2.2$ of \cite{P.-D. Liu1999}, we obtain that
 \begin{eqnarray}\label{gene-fiber entropy}
  h_{(\Phi\mathbf{m})}(f)=\frac{1}{k}\sum_{j=1}^{k}\int_{M} \log \text{Jac}(f_{j}(x))    d(\Phi\mathbf{m})(x).
 \end{eqnarray}
Then, by $(\ref{gene-skew measure entropy})$, $(\ref{gene-fiber})$ and $(\ref{gene-fiber entropy})$, we have
$$h_{(\Phi\mathbf{m})}(\mathcal{Q})=\log k + \frac{1}{k}\sum_{j=1}^{k}\int_{M} \log \text{Jac}(f_{j}(x))    d(\Phi\mathbf{m})(x).$$
It follows that
 $(\Phi\mathbf{m},\frac{1}{k}\sum_{j=1}^{k}\delta_{f_{j}(x)})$, $\forall x \in M$, is one of the equilibrium states for $P(T,\phi)$.

This completes the proof of Theorem \ref{thm:kernel-dense-function}.
\end{potB}

\section{The proof of Theorem C}

In this section, we study finitely-generated correspondences on torus and their equilibrium states with respect to real-valued bounded measurable potentials.

\begin{potC}

For simplicity, we will  prove the case where $m=1$. Suppose $f_{j}(x)=p_{j}x(\text{mod }1)$, $\forall x \in \mathbb{S}^{1}$, $1\leq j \leq k$.
The proof for the case where
$m \ge 2$ can be carried out similarly, and we leave the details to the readers.

 \emph{Step 1.}  Show that $P^{v}(T,\phi) \leq \log k$.

Let $\mathcal{Q}$ be an arbitrary transition probability kernel on $\mathbb{S}^{1}$ supported by $T$, and let
 $\mu$ be an arbitrary $\mathcal{Q}$-invariant probability measure on $\mathbb{S}^{1}$.
Note that $\phi$  defined in $(\ref{spe-poten})$ is an upper semi-continuous bounded simple function.
By the proofs of Theorem $9.10$ in \cite{Walters1982} and Lemma $2.11$ in \cite{Bowen1975}, we obtain the following inequality for $\phi$:

\begin{align*}\label{stil-hold}
 &\sup_{(\mu,\mathcal{Q})}\bigg\{h_{\mu}(\mathcal{\mathcal{Q}})+ \int_{X}\int_{T(x_{1})} \phi(x_{1},x_{2})d\mathcal{Q}_{x_{1}}(x_{2})d\mu(x_{1})\bigg\}
 \leq \lim_{\varepsilon \rightarrow 0}\limsup_{n \rightarrow +\infty}\frac{1}{n}\log \inf_{F} \notag   \\
 & \Big\{ \sum_{\underline{x}=(x_{1},x_{2},\cdots,x_{n+1})\in F}\exp{\big(\sum_{i=1}^{n}\phi(x_{i},x_{i+1})\big)}|F \text{ is a }(d_{n+1},\varepsilon)\text{-spanning set of }(\mathcal{O}_{n+1}(T),d_{n+1})\Big\}.
\end{align*}

 For any $\underline{x}=(x_{1},x_{2},\cdots,x_{n+1})\in \mathcal{O}_{n+1}(T)$, assume $x_{i+1} = f_{j_i}(x_i)$ with $1 \leq j_i \leq k$ and $1 \leq i \leq n$.
To show $P^{v}(T,\phi) \leq \log k$, it suffices to show
 \begin{align}\label{turn-only}
  \lim\limits_{\varepsilon\rightarrow 0} \limsup_{n\rightarrow +\infty}\frac{1}{n}\log \inf_{F} \Big\{\sum_{\underline{x}=(x_{1},x_{2},\cdots,x_{n+1})\in F}\frac{1}{\prod_{i=1}^{n}\text{Jac}(f_{j_{i}}(x_{i}))}\Big\} \leq \log k.
 \end{align}

 Since each $f_{j}(x)=p_{j}x(\text{mod }1)$, $1\leq j \leq k$, is expanding, then for any $x\in \mathbb{S}^{1}$, $x_{i} \in f_{j}^{-1}(x)$,
 \begin{equation}\label{pre num}
 \#f_{j}^{-1}(x)=p_{j}=\text{Jac}(f_{j}(x_{i})).
 \end{equation}
 Hence, we  take $0 < \varepsilon_{0} < \frac{1}{k}$ such that for any $0 < \varepsilon \leq \varepsilon_{0}$, $x \in \mathbb{S}^{1}$ and $1\leq j \leq k$,
 $$
 f_{j}^{-1}(B(x,\varepsilon))=\biguplus_{i=1}^{p_{j}}B(x_{i},\frac{\varepsilon}{p_{j}}),
 $$
 and the restriction $f_{j}|_{B(x_{i},\frac{\varepsilon}{p_{j}})}\colon B(x_{i},\frac{\varepsilon}{p_{j}})\to B(x,p_{j})$ is a homeomorphism.

 Let $K^{\prime}(\varepsilon)=\{x^{(l)} \in \mathbb{S}^{1}: 1 \leq l \leq K(\varepsilon)\}$ be an $\varepsilon$-net of $\mathbb{S}^{1}$.
We know that the set
 $$
 \mathcal{O}_{n+1}^{\prime}(T)=\{\underline{x}=(x_{1},\cdots,x_{n},x^{(l)}) \in \mathcal{O}_{n+1}(x^{(l)}): x^{(l)} \in K^{\prime}(\varepsilon)\}=\bigcup_{x^{(l)} \in K^{\prime}(\varepsilon)}\mathcal{O}_{n+1}(x^{(l)})
 $$is a $(d_{n+1},\varepsilon)$-spanning set of $(\mathcal{O}_{n+1}(T),d_{n+1})$.

 Hence, we have
 \begin{eqnarray*}
 \begin{split}
   &\lim\limits_{\varepsilon\rightarrow 0} \limsup_{n\rightarrow +\infty}\frac{1}{n}\log \inf_{F}\\
   & \Big\{ \sum_{\underline{x}=(x_{1},x_{2},\cdots,x_{n+1})\in F}\frac{1}{\prod_{i=1}^{n}\text{Jac}(f_{j_{i}}(x_{i}))}|F \text{ is a }(d_{n+1},\varepsilon)\text{-spanning set of }(\mathcal{O}_{n+1}(T),d_{n+1})\Big\}\\
   &\leq \lim\limits_{\varepsilon\rightarrow 0} \limsup_{n\rightarrow +\infty}\frac{1}{n}\log \{\max_{x^{(l)}\in K^{\prime}(\varepsilon)} \Phi_{n}(x^{(l)})\}.
 \end{split}
 \end{eqnarray*}

 Note that $\#\mathcal{O}_{n+1}^{\prime}(T) \leq K(\varepsilon)(\sum_{j=1}^{k}p_{j})^{n}$ and $\#\mathcal{O}_{n+1}(x^{(l)}) \leq (\sum_{j=1}^{k}p_{j})^{n}$. Combining these inequalities and $(\ref{pre num})$, we obtain that
 \begin{eqnarray*}
 \begin{split}
   &\lim\limits_{\varepsilon\rightarrow 0} \limsup_{n\rightarrow +\infty}\frac{1}{n}\log\big\{ \max_{x^{(l)}\in K^{\prime}(\varepsilon)} \Phi_{n}(x^{(l)})\big\}\\\leq &\lim\limits_{\varepsilon\rightarrow 0} \limsup_{n\rightarrow +\infty}\frac{1}{n}\log\bigg\{\sum_{j_{1},j_{2},\cdots,j_{k-1}=0}^{n}\binom{n}{j_{1}}\cdot \binom{n-j_{1}}{j_{2}}\cdot \cdots \cdot \binom{n-(j_{1}+j_{2}+\cdots+j_{k-2})}{j_{k-1}}\bigg\}\\
    =& \lim\limits_{\varepsilon\rightarrow 0} \limsup_{n\rightarrow +\infty}\frac{1}{n}\log k^{n}=\log k.
 \end{split}
 \end{eqnarray*}
Therefore, the inequality $(\ref{turn-only})$ holds, and hence $P^{v}(T,\phi) \leq \log k$.

\emph{Step 2.}  Show that  $P^{v}(T,\phi) \geq \log k$.

 Since $f_{j}(x)=p_{j}x(\text{mod }1)$,  the Lebesgue measure is  $f_{j}$-invariant for each $j$, $1 \leq j \leq k$. Let  $\mathcal{Q}_{x}=\frac{1}{k}\sum_{j=1}^{k}\delta_{f_{j}(x)}$, $\forall x \in \mathbb{S}^{1}$. Then the kernel $\mathcal{Q}$ is supported by $T$ and $\text{Leb}\mathcal{Q}=\text{Leb}$.
Hence, by $(\ref{variational topological pressure})$, we  obtain
 \begin{align}\label{int}
   P^{v}(T,\phi) \geq &P_{(\text{Leb},\frac{1}{k}\sum_{j=1}^{k}\delta_{f_{j}})}^{v}(T,\phi) \notag \\=& h_{\text{Leb}}(\mathcal{Q})+
  \int_{\mathbb{S}^{1}}\int_{T(x_{1})} \phi(x_{1},x_{2})d\big(\frac{1}{k}\sum_{j=1}^{k}\delta_{f_{j}(x_{1})}\big)(x_{2})d\text{Leb}(x_{1}) \notag   \\
  =& h_{\text{Leb}}(\mathcal{Q})+
  \int_{\mathbb{S}^{1} \setminus \{0\} }\int_{T(x_{1})} \phi(x_{1},x_{2})d\big(\frac{1}{k}\sum_{j=1}^{k}\delta_{f_{j}(x_{1})}\big)(x_{2})d\text{Leb}(x_{1}) \notag   \\
  =& h_{\text{Leb}}(\mathcal{Q})- \frac{1}{k}\log \prod_{j=1}^{k}p_{j}.
\end{align}

Similar to what we have done in the proof of Theorem \ref{thm:kernel-dense-function},
we have
\begin{eqnarray}\label{product project}
 h_{\text{Leb}}(\mathcal{Q})=h_{\widetilde{\pi}({\mathbb{P}_{k}}^{\mathbb{N}}\times \text{Leb})}(\sigma_{T})=h_{{\mathbb{P}_{k}}^{\mathbb{N}}\times \text{Leb}}(\widetilde{\sigma})=h_{{\mathbb{P}_{k}}^{\mathbb{N}}}(\sigma_{k})+h_{\text{Leb}}(f)=\log k+\frac{1}{k}\log \prod_{j=1}^{k}p_{j},
\end{eqnarray}where $\mathbb{P}_{k}=(\underbrace{\frac{1}{k}, \frac{1}{k}, \cdots, \frac{1}{k}}_{k})$ and $h_{\text{Leb}}(f)$ denotes the fiber-entropy of $\widetilde{\sigma}$ with respect to Lebesgue measure.

Then, by $(\ref{int})$ and $(\ref{product project})$, we have
$$
P^{v}(T,\phi) \geq P_{(\text{Leb},\frac{1}{k}\sum_{j=1}^{k}\delta_{f_{j}})}^{v}(T,\phi) = \log k.
$$
Combining step $1$ and step $2$, we have
$$
P^{v}(T,\phi)=P_{(\text{Leb},\frac{1}{k}\sum_{j=1}^{k}\delta_{f_{j}})}^{v}(T,\phi)= \log k.
$$
It follows that
 $(\text{Leb},\frac{1}{k}\sum_{j=1}^{k}\delta_{f_{j}(x)})$, $\forall x \in \mathbb{S}^{1}$, is one of the equilibrium states for $P^{v}(T,\phi)$.

This completes the proof of Theorem \ref{thmm-bounded-poten-equi}.
\end{potC}

\section{On the uniqueness of equilibrium states}

In this section, we discuss the uniqueness of equilibrium states via the  natural extensions. Let $T$ be a correspondence generated by finite $C^{2}$-expanding endomorphisms $f_{1},f_{2},\cdots, f_{k}$ on $M$.

\begin{remark}\label{unique-sense}
 In Theorem \ref{thm:kernel-dense-function}, we considered the existence for the correspondence $T$ and the locally H\"{o}lder continuous potential $\phi$. In the following, we will show that the equilibrium state $(\Phi\mathbf{m},\frac{1}{k}\sum_{j=1}^{k}\delta_{f_{j}})$ is unique in an appropriate sense. By $(\ref{topological pressure})$, $P(T,\phi)=P(\sigma_{T},\widetilde{\phi})$. Therefore, we  study the uniqueness of  the equilibrium state of correspondences using the left shift  map on the orbit space with respect to the corresponding potential $\widetilde{\phi}$,  given in (\ref{local constant}). Note that $\widetilde{\phi}\colon \mathcal{O}_{\omega}{(T)}\to \mathbb{R}$
   is also locally H\"{o}lder continuous (see Remark \ref{extend-holder}). We say a single-valued continuous  map $f\colon X\to X$ is \emph{topologically transitive} if for any non-empty open sets $U, V \subset X$, $\exists n \in \mathbb{N}$ such that $f^{n}(U) \cap V \neq \emptyset$.
 By Lemma \ref{open-transitive-expand}, we show that the correspondence $T$ is strongly transitive and distance-expanding. Hence,
$\sigma_{T}\colon \mathcal{O}_{\omega}{(T)}\to  \mathcal{O}_{\omega}{(T)}$ is topologically transitive and distance-expanding (see Proposition D.5, Proposition D.4 and Proposition D.7 of  \cite{Xiaoran Li-Zhiqiang Li-Yiwei Zhang2023}, respectively).  Therefore, applying  a  version  of  Ruelle-Perron-Frobenius Theorem (see \cite{Przytycki-Urbanski2010}, for example) to the  natural extension $(\mathcal{O}_{\omega}(T),\sigma_{T})$ and the potential $\widetilde{\phi}$, we conclude that $(\Phi\mathbf{m})(\frac{1}{k}\sum_{j=1}^{k}\delta_{f_{j}})^{\omega}$ is the unique equilibrium state for $P(\sigma_{T},\widetilde{\phi})$.
Let $(\mu, \mathcal{Q})$  be a equilibrium state for $P(T,\phi)$, where $\mathcal\mathcal{Q}_{x}=\sum_{j=1}^{k}P_{x,f_{j}}\delta_{f_{j}(x)}$, $\forall x \in M$.
Then $\mu(\sum_{j=1}^{k}P_{x,f_{j}}\delta_{f_{j}})^{\omega}$ is an equilibrium state for $P(\sigma_{T},\widetilde{\phi})$. By uniqueness, we have $$\mu(\sum_{j=1}^{k}P_{x,f_{j}}\delta_{f_{j}})^{\omega}=(\Phi\mathbf{m})(\frac{1}{k}\sum_{j=1}^{k}\delta_{f_{j}})^{\omega}.$$
This implies that $\mu=\Phi\mathbf{m}$. And
for $\Phi\mathbf{m}$-almost $x \in M$, $\forall A \in \mathcal{B}(M)$, we have $$\sum_{j=1}^{k}P_{x,f_{j}}\delta_{f_{j}(x)}(A)=\frac{1}{k}\sum_{j=1}^{k}\delta_{f_{j}(x)}(A).$$ This implies that, for $\Phi\mathbf{m}$-almost $x \in M$, $P_{x,f_{j}}=\frac{1}{k}$, $1 \leq j \leq k$. In this sense, we say $(\mu, \mathcal{Q})=(\Phi\mathbf{m},\frac{1}{k}\sum_{j=1}^{k}\delta_{f_{j}})$. Hence, we say $(\Phi\mathbf{m},\frac{1}{k}\sum_{j=1}^{k}\delta_{f_{j}})$ is the unique equilibrium state for $P(T,\phi)$.
\end{remark}

\begin{remark}\label{two-variation-coin}
  In Theorem \ref{thmm-bounded-poten-equi}, we considered the variational topological pressure $P^{\nu}(T,\phi)$ for $T$ and  $\phi \in B(\mathcal{O}_{2}(T),\mathbb{R})$. Note that if $\phi \in B(\mathcal{O}_{2}(T),\mathbb{R})$, then $\widetilde{\phi} \in B(\mathcal{O}_{\omega}(T),\mathbb{R})$.  By $(\ref{topological pressure})$, $P(T,\phi)=P(\sigma_{T},\widetilde{\phi})$, where $\phi\colon \mathcal{O}_{2}{(T)}\to \mathbb{R}$ and $\widetilde{\phi}\colon \mathcal{O}_{\omega}{(T)}\to \mathbb{R}$ are  continuous functions.
   A natural question arises:  Is there a similar notion of variational topological pressure $P^{v}(\sigma_{T}, \widetilde{\phi})$ for   the corresponding natural extension $(\mathcal{O}_{\omega}(T),\sigma_{T})$ and the potential  $\widetilde{\phi}$ such that $P^{v}(T, \phi) = P^{v}(\sigma_{T},\widetilde{\phi})$?  The answer is positive. Define  $P^{v}(\sigma_{T}, \widetilde{\phi})$ by
  $$
  P^{v}(\sigma_{T}, \widetilde{\phi}):=\sup_{\nu \in \mathcal{M}_{\sigma_{T}}(\mathcal{O}_{\omega}(T))}\bigg\{h_{\nu}(\sigma_{T})+\int_{\mathcal{O}_{\omega}(T)} \widetilde{\phi}(x_{1},x_{2},\cdots) d\nu\bigg\}.
  $$We will then  show that $P^{v}(T, \phi) = P^{v}(\sigma_{T},\widetilde{\phi})$.
From Lemma $A.9$ of \cite{Xiaoran Li-Zhiqiang Li-Yiwei Zhang2023}, for any $\phi \in B(\mathcal{O}_{2}(T),\mathbb{R})$, we have
$$
\int_{X}\int_{T(x_{1})} \phi(x_{1},x_{2})d\mathcal{Q}_{x_{1}}(x_{2})d\mu(x_{1})=\int_{\mathcal{O}_{\omega}(T)} \widetilde{\phi}(x_{1},x_{2},\cdots)    d\mu \mathcal{Q}^{\omega}(x_{1},x_{2},\cdots),
$$where $\widetilde{\phi}(x_{1},x_{2},\cdots)=\phi(x_{1},x_{2})$, $\mathcal{Q}$ is supported by $T$ and $\mu$ is $\mathcal{Q}$-invariant.
 Combining this with the fact that $h_{\mu}(\mathcal{Q})=h_{\mu \mathcal{Q}^{\omega}}(\sigma_{T})$, we obtain
$$
P^{v}(T, \phi) \leq P^{v}(\sigma_{T},\widetilde{\phi}).
$$
By  the proof of Proposition $6.18$ in \cite{Xiaoran Li-Zhiqiang Li-Yiwei Zhang2023},  if $\nu$ is an equilibrium state for $P^{v}(\sigma_{T},\widetilde{\phi})$, i.e.,
$h_{\nu}(\sigma_{T})+\int_{\mathcal{O}_{\omega}(T)} \widetilde{\phi}(x_{1},x_{2},\cdots) d\nu =P^{v}(\sigma_{T},\widetilde{\phi})$, then
there exists a transition probability kernel $\mathcal{Q}$ supported by $T$ and a $\mathcal{Q}$-invariant Borel probability measure $\mu$ such that $$h_{\nu}(\sigma_{T})+\int_{\mathcal{O}_{\omega}(T)} \widetilde{\phi}(x_{1},x_{2},\cdots) d\nu  \leq h_{\mu}(\mathcal{\mathcal{Q}})+ \int_{X}\int_{T(x_{1})} \phi(x_{1},x_{2})d\mathcal{Q}_{x_{1}}(x_{2})d\mu(x_{1}).$$
Thus, we conclude that
$P^{v}(T, \phi) = P^{v}(\sigma_{T},\widetilde{\phi})$.
\end{remark}

As  noted in the proof of Theorem \ref{thmm-bounded-poten-equi}, we have
  \begin{align*}
  &P^{v}(T, \phi)
  = h_{(\text{Leb}\mathcal{Q}^{\omega})}(\sigma_{T})+\int_{\mathcal{O}_{\omega}(T)} \widetilde{\phi}(x_{1},x_{2},\cdots)d(\text{Leb}\mathcal{Q}^{\omega})\\
  =&h_{\widetilde{\pi}({\mathbb{P}_{k}}^{\mathbb{N}}\times \text{Leb})}(\sigma_{T})+ \int_{\mathcal{O}_{\omega}(T)} \widetilde{\phi}(x_{1},x_{2},\cdots)d\widetilde{\pi}({\mathbb{P}_{k}}^{\mathbb{N}}\times \text{Leb}),
   \end{align*}where $\mathcal{Q}_{x}=\frac{1}{k}\sum_{j=1}^{k}\delta_{f_{j}(x)}, \forall x \in \mathbb{T}^{m}$. It follows that
  $\widetilde{\pi}({\mathbb{P}_{k}}^{\mathbb{N}}\times \text{Leb})$ is one of the equilibrium states for $P^{v}(\sigma_{T}, \widetilde{\phi})$.  We guess that such an equilibrium state is unique and hence the  equilibrium state $(\text{Leb},\frac{1}{k}\sum_{j=1}^{k}\delta_{f_{j}})$  is unique for $P^{v}(T, \phi)$.

\section*{Data availability statement}
No dataset is generated in this research.

\section*{Conflicts of interest statement}
The authors have no financial or non-financial conflicts of interest to declare that are relevant to this article.

\footnotesize{


\textit{ E-mail address:} \href{mailto:\Emaila}{zhangyu3@stu.xmu.edu.cn}

\textit{ E-mail address:} \href{mailto:\Emailb}{yjzhu@xmu.edu.cn}}

\end{document}